\documentclass[11pt,a4paper]{article}
 \usepackage{mathrsfs}
\usepackage[british]{babel}
\usepackage[utf8]{inputenc}
\usepackage{amsmath}
\usepackage{amssymb}
\usepackage{amsthm}
\usepackage{verbatim}
\numberwithin{equation}{section}
\usepackage{mathtools}
\usepackage{graphicx}
\usepackage[color=white]{todonotes}

\graphicspath{ {./images/} }

\newcommand{\R}{\mathbb R}
\newcommand{\E}{\mathbb E}
\renewcommand{\P}{\mathbb P}
\newcommand{\N}{\mathbb N}
\newcommand{\Z}{\mathbb Z}
\newcommand{\C}{\mathbb C}

\newcommand{\F}{{\cal {F}}}
\newcommand{\bF}{{\mathbb F}}
\newcommand{\<}{\left<}

 \renewcommand{\>}{\right>}
\newcommand{\abs}[1]{\left\vert#1\right\vert} 
\newcommand{\ind}{1\mkern-7mu1}

\newcommand{\Rp}{\R_+}
\newcommand{\B}{\mathscr B}
\newcommand{\SR}{{\mathscr S}(\R)}
\newcommand{\SprimeR}{{\mathscr S}'(\R)}
\newcommand{\supp}{\operatorname{supp}}
\newcommand\norm[1]{\left\lVert#1\right\rVert}

 \usepackage{hyperref}
\usepackage[framemethod=tikz]{mdframed}

\usepackage{amsthm}

\theoremstyle{plain}
\newtheorem{thm}{Theorem}[section]
\newtheorem{lem}[thm]{Lemma}
\newtheorem{cor}[thm]{Corollary}
\newtheorem{prop}[thm]{Proposition}

\newtheorem*{theorem*}{Theorem}
\theoremstyle{definition}

\newtheorem{notation}[thm]{Notation}
\newtheorem{assumption}[thm]{Assumption}
\newtheorem{ex}[thm]{Example}
\newtheorem*{asum*}{Assumption}

\newtheorem{rem}[thm]{Remark}

\renewcommand{\proof}{\vskip 0pt\noindent\textbf{\textit{Proof. }}}
\newcommand{\proofof}[1]{\vskip 0pt\noindent\textbf{\textit{Proof #1. }}}

\begin{document}
\title{\textbf{A generalized central limit theorem for critical marked Hawkes processes}}
\author{\Large{Anna Talarczyk}\\
University of Warsaw, Institute of Mathematics\\ ul. Banacha 2, 02-097 Warsaw, Poland\\ \hbox{e-mail:}
annatal@mimuw.edu.pl
}
\bigskip
\date{May 1, 2026}
\maketitle
\begin{abstract}
We prove a central limit type theorem for critical marked Hawkes processes. We study the case where the marks are i.i.d.\ with nonnegative values and their common distribution is either heavy tailed or has finite variance. The kernel function is of a multiplicative form and the mean number of future events triggered by a single event is $1$ (criticality).
We also assume that the base intensity function is heavy tailed. We prove convergence in law in the space of tempered distributions of the normalized empirical measure corresponding to the times of events. We also study convergence in law in the Skorokhod space of the normalized event counting process as the time is speeded up. 
In case when the distribution of marks is heavy tailed, the limit process is a stable process with dependent increments, while in case of finite variance, the limit process  is the same  Gaussian process as for the non marked Hawkes process.
We develop a new, robust method that may be applied to other self-exciting systems generalizing Hawkes processes.
For example we consider a non marked
self-exciting system where the number of excitations caused by single event is heavy tailed.
\end{abstract}

\bigskip

\\

 \textup{2020} 

\section{Introduction}\label{sec:intro}

Hawkes processes were first introduced by Hawkes in 1971 in \cite{Hawkes1} and \cite{Hawkes2}. The Hawkes process, is a type of a counting process  which is self-exciting,  in the sense that each arrival (jump) in the process  influences the intensity of future jumps. More precisely, let $N=(N_t)_{t\ge 0}$ denote the Hawkes process, then the values of  $N$ are nonnegative integers, $N$ is nondecreasing and for $\delta>0$ we have 
\begin{equation*}
 P(N_{t+\delta}-N_t=1|\F_t)=\Lambda(t)\delta+o(\delta),
\end{equation*}
where $\F_t$ denotes the sigma-field incorporating the whole history of the process up to time $t$ and $\Lambda(t)$ is the intensity, which itself depends on $N$, namely
\begin{equation*}
 \Lambda(t)=\mu(t)+\int_{(0,t)}\varphi(t-s)N(ds),
\end{equation*}
where $\varphi$ is a nonnegative function in $L^1(\R)$ such that $\varphi(s)=0$ for $s<0$ and $\mu$ is locally integrable. If we denote by $\tau_i$, $i=1,2,\ldots$ the succesive jump times of the Hawkes process, then $\Lambda$ can be written as
\begin{equation*}
 \Lambda(t)=\mu(t)+\sum_{0<\tau_i<t}\varphi(t-\tau_i).
\end{equation*}
The function $\mu$ is the intensity of  apperance of exogeneous events and $\varphi$ is responsible for the self exciting nature of the process. Often one takes $\mu(t)$ to be constant, as it will be the case in our model.

A natural extension of this process studied in the literature is the marked Hawkes process, where, additionally, each event (corresponding to the jump in the process) is assigned a certain mark and both the time of appearance of the event and its mark influence the intensity. The intensity takes the form
\begin{equation*}
 \Lambda(t)=\mu(t)+\sum_{0<\tau_i<t}\phi(t-\tau_i,\eta_i),
\end{equation*}
where, as before, $\tau_i$ is the time of $i$-th jump of the process $N$ and $\eta_i$ is the corresponding mark, $\phi$ is an appropriate kernel function. 
Very often it is assumed that the marks $\eta_i$, $i=1,2,\ldots$ are independent identically distributed and that each $\eta_i$ is independent of the sigma-field  generated by the process and the marks up to time $\tau_i$.
In our setup we will assume that the random variables $\eta_i$ take values in $\R_+$ and and are i.i.d. with common distribution $\nu$.
In this context one often considers  the function $\phi$ of the form
\begin{equation}
 \phi(t,x)=x\varphi(t)
 \label{e:kernel}
\end{equation}
for some fixed $\varphi$ and $\nu$ is a probability distribution on $\Rp$. This means that the events only differ by the importance/magnitude. 
Marked processes of  this form are used for example in 
epidemic type aftershock models (ETAS) in
seismology (see e.g. \cite{Ogata}) with functions $\varphi(t)=\frac 1{(1+t)^p}$. Such functions $\varphi$ will be studied later in our paper.

 The case $\nu=\delta_1$ corresponds to the usual (not marked) Hawkes process.

Due to their self-exciting nature, the Hawkes proceses and their generalizations, in particular marked Hawkes processes gained a lot of interest with many applications, providing natural models  e.g. in seismology, epidemiology, finance. There is a vast literature concerning Hawkes and marked Hawkes processes, see for example \cite{HawkesRev}, \cite{DaleyVereJones},  \cite{HorstXu}, \cite{HorstXu_marked}, \cite{Laub}, for and overview  of the results in literature.

Usually, the stable (subcritical) case  is considered, that is when $\int_{\R_+\times \R_+}\phi(s,x)ds\nu(dx)<1$. This implies that the average number of future excitations induced by a single event is smaller than $1$ and the intensity of the marked Hawkes process converges in law towards equilibrium. However there is an increasing interest in the nearly critical or critical cases \cite{Hardiman}, \cite{JaissonRosenbaum_nearly_unstable}, \cite{JR_2}, \cite{HorstXu}.

In the present paper we are interested in the central limit type result for the event counting process $(N_t)_{t\ge 0}$, that is, in the limiting behavior of 
\begin{equation}
 X_T(t)=\frac{N_{tT}-EN_{tT}}{F_T},\qquad t\ge 0
 \label{e:XTintro}
\end{equation}
as $T\to \infty$ and $F_T$ is an appropriate norming. 
This problem was studied first for subcritical Hawkes process by Hawkes and Oakes in \cite{HawkesOakes}. 
The authors showed that if $\int_0^\infty t\varphi(t)dt<\infty$, then, with the standard norming $F_T=\sqrt{T}$ the processes $X_T$ converge in law to a Brownian motion.
Bacry, Dellattre and Hoffmann in \cite{BacryDelattreHoffmann}  considered subcritical multivariate Hawkes processes, extending the result of \cite{HawkesOakes} and relaxing the condition on $\varphi$. Note that the multivariate Hawkes processes form a special case of marked Hawkes processes.
More general subcritical marked Hawkes processes and marked Hawkes random measures were investigated in 
\cite{HorstXu_marked}. See e.g. Section 1.1 of \cite{HorstXu} for a thorough and up to date literature review.

Much of the work was done for subcritical Hawkes processes, but there are also some recent results concerning critical and nearly critical cases. In \cite{JaissonRosenbaum_nearly_unstable} Jaisson and Rosenbaum considered 
 nearly critical  Hawkes processes. 
 They investigated the case when the kernel functions $\varphi_T$ depend on $T$ and are of the form $\varphi_T=a_T\varphi$ with $\norm{\varphi}_1=1$, $\int_{\Rp}t\varphi(t)dt<\infty$ and $a_T<1$, $a_T\to 1$. In the limit of \eqref{e:XTintro} they obtained the integrated Cox-Ingersoll-Ross process.
 The same authors in \cite{JR_2} investigated a similar, nearly critical case but with the kernel function having powerlaw tail: $\varphi(x)\sim x^{-1-\alpha}$, as $x\to \infty$, $\alpha\in (0,1)$. 
 They showed that the limits in law of $X_T$ are related to rough fractional diffusions.

 The recent work by Horst and Xu \cite{HorstXu} gives a detailed study of both subcritical and critical (not marked) Hawkes processes. The latter paper was an inspiration fot the present work. In particular, in \cite{HorstXu} the authors obtain a very interesting result in the case of the Hawkes process with heavy-tailed kernel function. In Theorem 2.14 of \cite{HorstXu} they show  that if  the kernel function $\varphi$ is such that 
\begin{equation}
\label{e:phi_cond}
 \int_t^\infty \varphi(s)ds\sim t^{-\alpha} \qquad as \ t\to \infty
 \end{equation}
 for some
 $\alpha\in(0,1)$ and $\int_{\R_+}\varphi(s)ds=1$
then, $\E N(T)$ behaves as $T^{1+\alpha}$ and  with the norming $F_T=T^{\frac{1+3\alpha} 2}$ the processes  \eqref{e:XTintro} converge in law in the Skorokhod space $D([0,\infty))$ to the process
\begin{equation}
 C\int_0^t (t-s)^{\alpha}s^{\frac \alpha 2}dB_s,
\label{e:HorstXu_lim}
 \end{equation}
where $B$ is a Brownian motion and $C$ is a constant. In fact Theorem 2.14 in \cite{HorstXu} is somewhat more general and allows the right hand side of \eqref{e:phi_cond} to be a more general function regularly varying at infinity.

Particularly interesting about this result is that the fact that $\varphi$ is heavy tailed implies that the limit process exhibits long range dependence.

 It is also worth observing that the covariance function of the Gaussian process in \eqref{e:HorstXu_lim} is, up to a constant, of the form
\begin{equation*}
 K(s,t)=\int_0^\infty G^{(\alpha)}\left((G^{(\alpha)} \ind_{[0,s]})(G^{(\alpha)}\ind_{[0,t]})\right)(x)dx
 =C\int_0^\infty\left((G^{(\alpha)} \ind_{[0,s]})(x)(G^{(\alpha)}\ind_{[0,t]})(x)\right)x^\alpha dx,
\end{equation*}
where $G^{(\alpha)}$ is the potential of the $\alpha$-stable L\'evy process totally skewed to the right. This bears some resemblance to the limits obtained for occupation time fluctuations of $\alpha$-stable branching particle systems e.g. in \cite{functlim4}.

The motivation for the present work was to investigate how the long range dependence phenomenon translates to the case of marked Hawkes processes. In particular the case where the marks are heavy tailed, as well as understanding better the relation with  branching particle systems. 

In this paper we study marked Hawkes processes with  kernel functions of the form
\eqref{e:kernel}.
We assume that  $\norm{\varphi}_1=\int_{\R_+}\varphi(x)dx=1$, $\varphi$ satisies \eqref{e:phi_cond} and the marks are i.i.d. with the common distribution $\nu$ with $\int_{\Rp} x\nu(dx)=1$. We show that if $\nu$ is heavy-tailed and is in the normal domain of attraction of an $(1+\beta)$-stable law with $1>\beta>\alpha$, then the suitable norming is 
\begin{equation*}
 F_T=T^{\frac{ 1+\alpha(2+\beta)}{1+\beta}}
\end{equation*}
and the processes $X_T$ defined in \eqref{e:XTintro} converge in the sense of finite dimensional distributions to a $(1+\beta)$-stable process of the form
\begin{equation}
 C\int_{[0,t]} (t-s)^{\alpha}s^{\frac \alpha {1+\beta}}dL_{1+\beta}(s),
 \label{e:limit_proc_2}
\end{equation}
where $L$ is a $(1+\beta)$-L\'evy process totally skewed to the right (see Theorem \ref{thm:main}b)). The limit process \eqref{e:limit_proc_2} has a continuous modification (see Proposition \ref{prop:continuity}). Under an additional, technical, but quite natural condition we prove convergence in law in the Skorokhod space $D([0,\infty))$ of c\`adl\`ag functions equipped with $J_1$ topology (Theorem
\ref{thm:tightness}). 

In Theorem \ref{thm:main}a) we show
 convergence of normalized empirical measures
\begin{equation*}
 X_T=\frac 1{F_T}\left(\sum_{k}\delta_{\frac{\tau_k} T}
 - E \sum_{k}\delta_{\frac{\tau_k} T}\right)
\end{equation*}
in law in $\SprimeR$  -- the space of tempered distributions on $\R$. Here, as before, $\tau_k$ denote the jump times of the process $N$. Note that in this notation $\<X_T,\ind_{[0,t]}\>=X_T(t)$.

We show that $X_T$ converges in law in $\SprimeR$ to a $(1+\beta)$-stable random variable $X$ such that
for any $f\in \SR$ the real valued random variable $\<X,f\>$
is $(1+\beta)$-stable
and it has the same law as 
\begin{equation*}
 C\int_0^\infty G^{(\alpha)}f(x)x^{\frac \alpha{1+\beta}} L_{1+\beta}(dx),
\end{equation*}
where $C$ is a positive constant, that does not depend on $f$, $G^{(\alpha)}$ denotes the potential of an $\alpha$-stable L\'evy process totally skewed to the right and $L_{1+\beta}$ is a $(1+\beta)$-stable random measure with Lebesgue control measure and skewness intensity $1$ (cf. \cite{ST} for the definitions). 

In the case when $\nu$ has finite variance the limit process for $(X_T(t))_{t\ge 0}$ is the same as for the not marked Hawkes process in Theorem 2.14 of \cite{HorstXu}, that is, it has the form \eqref{e:HorstXu_lim}. We also have convergence of $X_T$ in law in $\SprimeR$ to a centered Gaussian random variable $X$, such that for $f\in \SprimeR$ the variance of $\<X,f\>$ has the form
\begin{equation*}
 C\int_0^\infty (G^{\alpha}f(x))^2x^\alpha dx.
\end{equation*}

It should be stressed that our approach is different from that of Horst and Xu \cite{HorstXu}. The method of \cite{HorstXu} consisted in the study of the solutions to equations with respect to a Poisson random measure and using convergence theorems for martingales. We rely very strongly on the branching property of the marked Hawkes process.
We prove convergence by working with Laplace transforms and making use of the branching representation of marked Hawkes processes. 
A major difficulty is also caused by the fact that we consider marked Hawkes processes and that we allow the marks to be heavy tailed. Moreover, our method may be extended beyond the world of Hawkes processes and may be applied to other self-exciting processes (see e.g. Theorem \ref{thm:beta_branching}).

Practically from the beginning, when the Hawkes processes were introduced it was clear (see \cite{HawkesOakes}), that these processes may be equivalently represented by a branching system. See \cite{DaleyVereJones} Ex. 6.4(c) or \cite{Li_Pang} for the branching description of the marked Hawkes process. 
Since each event increases the intensity, one can think of all the subsequent events induced by this part of intensity as the ``offspring'' of that event.
In our case, 
since the kernel function of the form \eqref{e:kernel},  the marks are i.i.d. and $\norm{\phi}_1=1$ we have some exogeneous events which are distributed according to a Poisson random measure with intensity $\mu \ell$, where $\ell$ denotes Lebesgue measure. Then an event which happens at time $t_i$ is  assigned independently a mark $\eta_i$. 
This event then has a random number $\theta_i$ of ``children'' which, conditionally on $\eta_i=x$ is Poisson with parameter $x\int_{\Rp}\varphi (s)ds$, the children are placed at points $t_i+\xi_{i,k}$, and are assigned i.i.d. marks $\eta_{i,k}$, respectively, $k=1,\ldots, \theta_i$ where $(\eta_{i,k})_k$ are i.i.d. with law $\nu$ and
$(\xi_{i,k})_k$ are i.i.d with density $\frac 1{\norm{\varphi}_1}\varphi$, both sequences are independent of each other.
This branching model is described in detail in the next section.

It is also worth mentioning that the limiting process and the random variable $X$ resembles the spatial structure of limits that we  obtained together with Bojdecki and Gorostiza in \cite{functlim4} and subsequent study of   the occupation time fluctuations of a  branching particle system where the particles moved according to an $\alpha$-stable L\'evy process. Our result was a generalization of an earlier result of Dawson Gorostiza and Wakolbinger in \cite{DGW}.  We comment on the relation of those branching particle system models and the branching system related to the marked Hawkes process in Remark \ref{rem:BGT}.

In Theorem \ref{thm:beta_branching} we show that the limiting behavior of the number of events in the marked Hawkes process with heavy tailed marks is the same as that of a similar branching system without marks but where the distribution  number of offspring is heavy tailed.

The paper is organized as follows:
In Section \ref{sec:result} we first describe in detail the branching representation of the marked Hawkes process and we present our main results.
In Section \ref{sec:Lapl_tr_and_mean} we derive the Laplace transform for $X_T$. Section \ref{sec:proofs} contain the proofs of the main results.

\section{Results}\label{sec:result}

As mentioned in the Introduction, we will use the interpretation of a marked Hawkes processes in terms of a branching particle system on $\Rp\times \Rp$. Each event in 
the marked Hawkes process occurring at time $t$ and having a mark $x$ is represented by a ``particle'' at $(t,x)$.

The particle system is determined by several parameters which we describe presently.

\medskip
\textbf{Basic assumptions and notation.} \label{basic_assumptions}
 Let $\mu>0$ be the basic (exogenous) intensity which is constant in time.
 In the present paper we will only consider the \textbf{critical case} $\norm{\phi}_1=1$. Moreover, we consider the kernel function of a particular form:
\begin{equation}
 \phi(t,x)=x\varphi(t),\label{e:form_of_phi}
\end{equation}
where $\varphi:\R\mapsto \R_+$ is a Borel function vanishing on $(-\infty,0)$ and such that 
\begin{equation}
\label{e:critical}
\int_0^\infty \varphi (t)dt=1.
\end{equation}
 Moreover,  we assume that the probability measure $\nu$ has mean $1$:
\begin{equation}
 \int_0^\infty x\nu(dx)=1.
 \label{e:Enu}
\end{equation}
This implies that $\norm{\phi}_1=1$.

By $\xi_1,\xi_2,\ldots$ we denote a sequence of i.i.d. random variables with density $\varphi$ and  by $\eta_1, \eta_2, \ldots$ i.i.d. random variables with distribution $\nu$.  Recall that $\E \nu=1$ by \eqref{e:Enu}. The sequences $(\xi_j)_j$ and $(\eta_j)_j$ are independent.
We will also use the notation $\xi$ and $\eta$  for generic independent random variables of this type.

\medskip
Positive constants, are denoted by $C,C_1,C_2,\ldots$ and may change from line to line.

\medskip

The particle system corresponding to a marked Hawkes process is constructed in the following way:

Suppose that the exogenous ($0$-generation) particles are distributed according to a Poisson random measure on $\Rp\times \Rp$ with intensity $\mu\ell\otimes \nu$.
(Note that this is equivalent to taking a Poisson random measure on $\Rp$ with intensity $\mu\ell$ where $\ell$ denotes the Lebesgue measure and assigning to each point of this random measure independently a ``mark'' from distribution $\nu$).

Subsequently, each of the particles independently gives rise to some  offspring. A particle present at $(t,x)$ produces a random number of offspring in the next generation. The number of immediate offspring is described by a random variable $\theta_x$ which is Poisson with parameter $x\int_0^\infty \varphi(s)ds=x$. The children particles are displaced with respect to the parent in the following way: given the parent at $t$, the offspring particles are at 
$t+\xi_j$, $j=1,\ldots \theta_x$, where $\xi_j,$ $j=1,2, \ldots$ are i.i.d. with density $\varphi$ with the respect to the Lebesgue measure, and they are independent of  $\theta_x$. The children particles then branch according to the same mechanism. Each particle branches independently of all the other particles in the system.

We denote by $M$ the empirical measure of the system of all the exogenous particles and all their  offspring of all generations, that is, for $A\in \B(\Rp\times\Rp)$ 
\begin{equation*}
 M(A)=\{\# \textrm{of particles of all generations in the set $A$}\}.
\end{equation*}
For $B\in \B(\R_+)$ let us also denote 
\begin{equation*}
  N(B)=M(B\times \R_+).
\end{equation*}
The points of the random point measure $N$ are the times 
of occurrence of events in the associated marked Hawkes process.

In this setup the mean number of immediate offspring of a single particle is\\ $\int_0^\infty\int_0^\infty x\varphi(s)\nu(dx)ds=1$. Therefore the branching system is critical. 
 This implies that the total number of offspring of all generations of a single particle is a.s. finite. As there is only a finite number of exogeneous particles in $[0,t]$, this shows that $ N ([0,t])<\infty$ a.s.

Note that in this way the counting process $(N(t))_{t\ge 0}$ defined as 
\begin{equation*}
 N(t):=N([0,t]), \qquad t\ge 0
\end{equation*}
is quite a natural marked Hawkes process in which  an event (a jump) at time $t$ is independently assigned a mark from distribution  $\nu$ and if the mark is $x$ then this event increases the intensity of future events (jumps) by $x\varphi(s-t)$, $s\ge 0$, starting from $N_0=0$.
Also observe, that if $\nu=\delta_1$, then $(N(t))_{t\ge 0}$ is the usual unmarked critical Hawkes process.

The random point measure $N$ may be considered on the whole real line, although $N((-\infty,0))=0$.
For a Borel function $f:\R\mapsto \R$ we write
\begin{equation*}
 \<N,f\>=\int_0^\infty f(t) N(dt)
\end{equation*}
whenever the integral makes sense. In particular, we can always write it when $f$ is nonnegative, although the integral may be infinite. As observed earlier, it is always well defined if $f$ has a bounded support.

\medskip

 Apart from the basic assumptions \eqref{e:form_of_phi} -- \eqref{e:Enu} in our main theorem we will impose the following assumptions on $\varphi$ and $\nu$:
 
 Let us denote 
 \begin{equation*}
  \Phi(t)=\int_t^\infty \varphi(s)ds.
 \end{equation*}

\begin{assumption} \label{ass:phi_alpha}
The Borel function $\varphi:\Rp\mapsto \Rp$ satisfies \eqref{e:critical} and exists $\alpha\in(0,1)$
and $c_\varphi>0$ such that 
\begin{equation*}
 \lim_{t\to \infty } t^{\alpha} \Phi(t)=c_\varphi.
\end{equation*}
\end{assumption}

\begin{assumption} $\nu$ is a probability measure satisfying \eqref{e:Enu} and either
\begin{itemize}
 \item [($\beta$)]
 \label{ass:nu_beta}
 there exist $\beta\in(0,1)$ and $c_\nu>0$ such that
 \begin{equation}
  \lim_{x\to\infty}x^{1+\beta}\int_x^\infty \nu(dz)=c_\nu
\label{e:nu_beta}
  \end{equation}
\item[or]
\item[($2$)]
\begin{equation*}
 \int_{\Rp} x^2\nu(dx)<\infty.
\end{equation*}
\end{itemize}
\end{assumption}

\begin{rem}
 Assumption \ref{ass:phi_alpha} implies that the probability distribution with the density $\varphi$ is in the normal domain of attraction of an $\alpha$-stable distribution totally skewed to the right (cf. \cite{Feller}, Ch. XVII.5 Theorem 2). That is, if $\xi_1,\xi_2, \ldots$ are i.i.d. with density $\varphi$ then
 \begin{equation*}
  \frac 1{n^{1/\alpha}}\sum_{k=1}^n \xi_k\Rightarrow CL_\alpha(1) \end{equation*}
as $n\to \infty$, where $C>0$ is a constant, and $L_\alpha(1)$ has an $\alpha$-stable distribution with Laplace transform of the form
\begin{equation}
 Ee^{-\lambda L_\alpha(1)}=e^{-{\lambda^\alpha}/{\cos(\frac {\pi\alpha}2)}}\qquad \lambda\ge 0
 \label{e:Lapl_alpha_stab}
\end{equation}
and characteristic function
\begin{equation}
 E{e^{i\theta L_\alpha(1)}}=\exp\left\{-|\theta|^\alpha\left(1-i({\rm sgn}\, \theta)\tan\frac{\pi}{2}\alpha\right)\right\}
 \label{e:char_alpha_stable}
\end{equation}
(b) Assumption \ref{ass:nu_beta}(b) implies that $\nu$ is in the normal domain of attraction of the $(1+\beta)$-stable law, totally skewed to the right, that is 
\begin{equation*}
  \frac 1{n^{1/(1+\beta)}}\sum_{k=1}^n(\eta_k-1) \Rightarrow CL_{1+\beta}(1), 
\end{equation*}
where $C>0$ is a constant and $L_{1+\beta}(1)$ 
Laplace transform and characteristic functions described by the same formulas \eqref{e:Lapl_alpha_stab} and \eqref{e:char_alpha_stable}, but with $\alpha$ replaced by $(1+\beta)$. Note that now $\cos(\frac{\pi(1+\beta)}2)$ is negative.
% 
% has Laplace transform
% \begin{equation}
% \label{e:Lapl_beta}
%  Ee^{-\lambda L_{1+\beta}(1)}=e^{-\lambda^{1+\beta}/\cos(\frac{\pi(1+\beta)}{2})}\qquad \lambda\ge 0
% \end{equation}
% (note that now $\cos(\frac{\pi(1+\beta)}2)<0$) and characteristic function
% \begin{equation*}
%  E{e^{i\theta L_{1+\beta}(1)}}=\exp\left\{-|\theta|^{1+\beta}\left(1-i({\rm sgn}\, \theta)\tan\frac{\pi}{2}(1+\beta)\right)\right\}
% \end{equation*}
\end{rem}

\textbf{Stable L\'evy processes.}

Later on we will also use the notation $(L_\alpha(t))_{t\ge 0}$ and $(L_{1+\beta}(t))_{t\ge 0}$ for c\`adl\`ag stable L\'evy processes (processes with stationary independent increments) with $L_\alpha(1)$ and $L_{1+\beta}(1)$  as above. Note that these are spectrally positive processes - they do not have negative jumps. If $\alpha\in(0,1)$ then the process $(L_\alpha(t))_{t\ge 0}$ is nonnegative and nondecreasing, however $(L_{1+\beta}(t))_{t\ge 0}$ does not have these properties. In fact, the a.e. sample paths of  $L_{1+\beta}$ have infinite variation.

Our main references on stable random variables, stable processes and stable random measures are \cite{ST} and \cite{Feller}.

\medskip

\textbf{The resolvent.}

\nopagebreak
We denote 
\begin{equation}
 R(s)=\sum_{k=1}^\infty \varphi^{(*k)}(s),
\label{e:R}
 \end{equation}
where $\varphi^{(*k)}$ denotes the $k$-th convolution. $R$ is called the resolvent and 
 plays an important role in the context of Hawkes processes. This is also the case in the present paper.  

Similarly as in \cite{HorstXu} we denote
\begin{equation}
 \label{e:IR_def}
 I_R(t)=\int_0^t R(s)ds.
\end{equation}

From (8.6.3) in \cite{Bingham} (see also \cite{HorstXu}, Proposition 2.5 (3)) 
it follows that under  Assumption \ref{ass:phi_alpha} we have
\begin{equation}
 \label{e:IR}
\lim_{T\to \infty}\frac 1{T^\alpha} I_R(T)=c_\alpha,
\end{equation}
where
\begin{equation}
 \label{e:c_alpha}
 c_\alpha= \frac 1{c_\varphi\Gamma(1+\alpha)\Gamma(1-\alpha)}.
\end{equation}

\textbf{Normalized random field.}

Let $T\ge 1$ and assume $F_T$ is a (deterministic) positive norming depending on $T$, to be chosen later.
For any Borel function $f:\R\mapsto \R$ we denote
\begin{equation}
\label{e:fT}
 f_T(t)=\frac 1{F_T}f(\frac t T), \qquad t\ge 0.
\end{equation}

We are interested in the behavior as $T\to \infty$ of random fields $N_T$ defined in the following way:
\begin{equation}
\label{e:NT}\<N_T,f\>:=\<N,f_T\>,
\end{equation}
for any nonnegative Borel function $f$.

\begin{rem}(a)
The behavior of $N_T$ gives information on behavior of the 
process counting the events of the marked Hawkes process
 with time scaled by $T$. For any  $n\in \N$, $\theta_1,\ldots, \theta_n\in \R$ and $t_1,\ldots, t_k\in \R_+$
if $f$ is of the form
\begin{equation}
f(t)=    \sum_{k=1}^n \theta_k\ind_{[0, t_k]}
\label{e:f_sum}
\end{equation}
then
 \begin{equation*}
  \<N_T,f\>=\frac 1{F_T}\sum_k \theta_k N(t_kT),
 \end{equation*}
 \end{rem}

Our aim is to investigate convergence in law of
\begin{equation}
 X_T:=N_T-EN_T
\label{e:XT}
 \end{equation}
in an appropriate sense as $T\to \infty$.

For suitably chosen norming we show convergence of $X_T$ in law in the space of tempered distributions $\SprimeR$, dual to the Schwartz space $\SR$ of smooth and quickly decreasing functions. In this context $\<\cdot, \cdot\>$ denotes the duality between $\SprimeR$ and $\SR$.

One of the main ingredients needed to show this convergence is the  study of convergence in law of $\<X_T,f\>$ for $f\in \SR$. However, we will show the convergence of $\<X_T,f\>$ for a larger class of functions.  In particular, taking functions $f$ of the form \eqref{e:f_sum} we obtain convergence of finite dimensional distributions of the normalized marked Hawkes process
\begin{equation}
 X_T(t)=\frac 1{F_T} \left(N(Tt)-E N(Tt)\right) \qquad t\ge 0.
\label{e:XTproc}
 \end{equation}

 To study at the same time both the convergence of $X_T$ in law in $\SprimeR$ and the convergence of finite dimensional distributions of the process $(X_T(t))_{t\ge 0}$  we introduce the following classes of functions:

\begin{notation}\label{defn:Fgamma}
Let $\gamma >0$. By  $\bF_\gamma$ we denote the class of functions  $f:\Rp\mapsto \R$ which have at most finite number of points of  discontinuity and 
such  that there exists  $C>0$ (which may depend on $f$ and $\gamma$) such that 
\begin{equation}
\label{e:Fgamma}
 \abs{f(s)}\le \frac C{(1+s)^\gamma}, \qquad \forall s\in \Rp.
\end{equation}
By $\bF_\gamma^+$ we will denote the set of nonegative functions belonging to $\bF_\gamma$.
As usual, whenever needed,  we extend $f$ to $\R$ by setting $f(t)=0$ for $t<0$. 
Let us also denote 
\begin{equation}
 \bF =\bigcap_{\gamma>0}\bF_\gamma \qquad \textrm{and}\qquad \bF^+ =\bigcap_{\gamma>0}\bF_\gamma^+.
\end{equation}
\end{notation}
Observe that the functions of the form \eqref{e:f_sum} belong to $\bF$. Also, if $f\in \SR$ (the Schwartz space), then $f$ restricted to $\Rp$ belongs to $\bF$.

For functions $f\in \bF$ the expectation $E\<N,f\>$ is well defined and finite. More precisely, we have the following:
\begin{prop}\label{prop:expectation}
Suppose that  Assumption \ref{ass:phi_alpha} is satisfied and that $f\in \bF_{\gamma}$ for some $\gamma>1+\alpha$. Then
\begin{equation}
\label{e:NF}
 E\<N,|f|\><\infty.
\end{equation}
Moreover, 
\begin{align}
 \lim_{T\to \infty}\frac {1}{T^{1+\alpha}} E\<N, f(\frac \cdot T)\>
=&\alpha c_\alpha\int_0^\infty G^{(\alpha)}f(t) dt
\label{e:meanGalpha}\\
=&c_\alpha\int_0^\infty t^\alpha f(t)dt,
\label{e:talphaF}
 \end{align}
 where 
 \begin{equation}G^{(\alpha)}f(t)=\int_0^\infty f(t+s)s^{\alpha-1}ds
  \label{e:Galpha}
 \end{equation}
 and $c_\alpha$ is as in \eqref{e:c_alpha}.
\end{prop}

Although \eqref{e:talphaF} is simpler, the term on the right hand side of \eqref{e:meanGalpha} perhaps better explains why the limit is of this form. Namely, observe that up to a constant $G^{(\alpha)}$ is the $\alpha$-stable potential related to an $\alpha$-stable L\'evy process totally skewed to the right.

Suppose that  $(L_\alpha(t))_{t\ge 0}$ is an $\alpha$-stable L\'evy process totally skewed to the right with $L_\alpha(1)$ having the Laplace transform of the form \eqref{e:Lapl_alpha_stab}. 
Let $p_t$ denote the density of $L_\alpha(t)$. It is well known that $p_t$ is supported on $(0,\infty)$ and furthermore, due to the self-similarity property of the process $(L_\alpha(t))_{t\ge 0}$ for any $a>0$ we have
\begin{equation*}
 p_{at}(x)=a^{-\frac 1\alpha}p_t(xa^{-\frac 1\alpha}).
\end{equation*}
We also have
\begin{equation*}
 p_1(x)\le \frac C{1+\abs{x}^{1+\alpha}}.
\end{equation*}

Using these two properties it is straightforward to check  that if $\alpha\in(0,1)$, then  for any $y>0$ we have
\begin{equation*}
 \int_0^\infty p_s(y)ds=C y^{\alpha-1}
 \end{equation*}
for some constant $C>0$. Hence,
\begin{equation*}
 \int_0^\infty Ef(t+L_\alpha(s))ds=CG^{(\alpha)}f(t)
\end{equation*}

The reason why the stable potential appears in \eqref{e:meanGalpha}
is that if $\xi_1,\xi_2,\ldots$ are i.i.d with density $\varphi$ then the Assumption \ref{ass:phi_alpha} implies that the processes 
\begin{equation*}
 \frac{\xi_1+\ldots +\xi_{\lfloor tT^\alpha\rfloor}}{T}, \qquad t\ge 0
\end{equation*}
converge in law in $D([0,\infty))$ equipped with $J_1$ topology to the process $CL_\alpha(t))_{t\ge 0}$ as $T\to \infty$ (cf. \cite{Feller} Ch. XVII.5, Theorem 2 on p. 577, and \cite{Skorokhod57}, Theorem 2.7). 
Although the proof of Proposition  \ref{prop:expectation} does not refer directly to this fact, in Remark \ref{rem:heuristics} we give a  short heuristic argument explaining the role of the potential
 $G^{(\alpha)}$  in this setup.

\bigskip
We now proceed to our main result. Theorem \ref{thm:main} below
gives convergence of the normalized $\SprimeR$ random variables $X_T$, as well as convergence of finite dimensional distributions of the process $(X_T(t))_{t\ge 0}$.
The functional convergence in the Skorokhod space $D([0,\infty))$ discussed in Theorem \ref{thm:tightness}.

\begin{thm}
 \label{thm:main}
Let $\varphi$ and $\nu$ satisfy the Assumptions \ref{ass:phi_alpha} and \ref{ass:nu_beta}($\beta$) with $\alpha<\beta$ and let $F_T$ be the norming given by
\begin{equation}
 \label{e:FT} 
 F_T=T^{\frac{1+\alpha(2+\beta)}{1+\beta}}.
\end{equation}
Then \\
a) $X_T$ defined by \eqref{e:XT}
converge in law in $\SprimeR$ as $T\to \infty$ to an $\SprimeR$ valued $(1+\beta)$-stable random variable $X$ determined by
\begin{multline}
\label{e:f_char_X}
E{\exp}\{i\langle X,f\rangle\}\\
={\rm exp}\biggl\{-\mu K(-\cos \frac {\pi (1+\beta)}2)\int_{\Rp}|G^{(\alpha)}f(x)|^{1+\beta}\biggl(
1-i({\rm sgn}G^{(\alpha)}f(x))\tan\frac{\pi}{2}(1+\beta)\biggr)x^{\alpha}dx\biggr\}, \\
\qquad f\in \SR
\end{multline}
where
$$K=\frac {c_\nu}{\beta}\Gamma(1-\beta)c_\alpha^{2+\beta}\alpha^{1+\beta}$$
and $c_\alpha$ is given by \eqref{e:c_alpha}.
\\
b) The processes $(X_T(t))_{t\ge 0}$ defined by \eqref{e:XTproc} converge in the sense of finite dimensional distributions as $T\to \infty$ to the process $\frac {(\mu K(-\cos \frac {\pi (1+\beta)}2))^{ 1/{1+\beta}}}\alpha \zeta$, where $\zeta$ is given by
\begin{equation}
\label{e:limit_proc}
\zeta_t= \int_{[0,t]} u^{\frac \alpha{1+\beta}}(t-u)^\alpha L_{1+\beta}(du).
\end{equation}
Here $(L_{1+\beta}(t))_{t\ge 0}$ stands for a $(1+\beta)$-stable L\'evy process totally skewed to the right.
\\
\end{thm}

\begin{rem} \label{rem:main} a) Let us observe that the right hand side of \eqref{e:f_char_X} is a characteristic function of an integral with respect to a $(1+\beta)$-stable random measure. The definition and the properties of such integrals are discussed in detail in Chapter 3 of \cite{ST}.  

More precisely, for any fixed $f\in \SR$ the real-valued random variable $\<X,f\>$ has the same distribution as
\begin{equation*}
 \left(K\mu(-\cos \frac {\pi (1+\beta)}2)\right)^{1/(1+\beta)} \int_0^\infty G^{(\alpha)}f(u)u^{\frac \alpha{1+\beta}}L_{1+\beta}(du),
\end{equation*}
 where, in the terminology of \cite{ST}, $L_{1+\beta}(du)$ denotes the integration with respect to $(1+\beta)$-stable random measure with Lebesgue control measure and skewness intensity $1$ (which in fact is the same as integrating with respect to the L\'evy process $L_{1+\beta}$ as in part b)).  

This integral is well defined, since  if  $f\in \bF_\gamma$ with $\gamma>\frac{2+\beta}{1+\beta}(1+\alpha)$ then  
\begin{equation*}
 \int_0^\infty\abs{G^{(\alpha)}f(t)}^{1+\beta}t^{\alpha}dy<\infty,
\end{equation*}
and clearly $\SR\subset \bF_\gamma$ for any $\gamma>0$.

Therefore, for fixed $f\in \SR$,  $\<X,f\>$ has
 $(1+\beta)$-stable law with  scaling parameter $$\sigma_f=\left(\int_0^\infty  K\mu(-\cos \frac {\pi (1+\beta)}2) \abs{G^{(\alpha)}f(t)}^{1+\beta}t^\alpha dt\right)^{\frac 1{1+\beta}},$$
skewness  parameter equal to $$
\frac 1{\sigma_f^{1+\beta}}\left(-K\mu \cos \frac {\pi (1+\beta)}2 \int_0^\infty\abs{G^{(\alpha)}f(t)}^{1+\beta}
\operatorname{sgn}(G^{(\alpha)}f(t)) t^\alpha dt\right)
$$
and shift parameter equal to $0$
(cf. Property 3.2.2 in \cite{ST}).

Let us also observe, that if $f$ is nonnegative, then $\<X,f\>$ is totally skewed to the right and therefore its Laplace transform is finite and determines the distribution, moreover
\begin{equation}
 E\exp\left\{-\<X,f\>\right\}=\exp\left\{
  K\mu\int_0^\infty \left(G^{(\alpha)}f(t)\right)^{1+\beta} t^\alpha dt
  \right\}
  \label{e:Lapl_X}.
\end{equation}
See Proposition 1.2.12 in \cite{ST}.

b) We chose to write the process $\zeta$ in the form \eqref{e:limit_proc} in order to show the analogy to the result of \cite{HorstXu}.  If we formally set $\beta=1$ and replace $L_{1+\beta}$ by the Brownian motion then $\zeta$ the formula \eqref{e:limit_proc} gives the limit process obtained by Horst and Xu (\cite{HorstXu}, Thm. 2.14).

Our process $\zeta$ may be also expressed in terms of $G^{(\alpha)}$. It is trivial to observe that 
\begin{equation}
 G^{(\alpha)}\ind_{[0,t]}(u)=\int_0^\infty \ind_{[0,t]}(u+y)y^{\alpha-1}dy =\int_0^{t-u}y^{\alpha-1}dy\ind_{[0,t]}(u)=\frac 1\alpha (t-u)^\alpha\ind_{[0,t]}(u).
 \label{e:5.4a}
\end{equation}
Hence the stable process $\zeta$ defined by \eqref{e:limit_proc} can be written as
\begin{equation}
 \zeta_t=\alpha \int_{[0,t]} G^{(\alpha)}(\ind_{[0,t]})(u)
 u^{\frac \alpha {1+\beta}}L_{1+\beta}(du),\label{e:5.3}
\end{equation}
therefore parts a) and b) of Theorem \ref{thm:main} are consistent. 

\medskip
c)
Regarding the condition $\alpha<\beta$ i  Theorem \ref{thm:main}, observe that if we formally set $\alpha=\beta$ in \eqref{e:FT} then $F_T=T^{1+\alpha}=T^{1+\beta}$, which is the same norming as for the mean (cf. Proposition \ref{prop:expectation}). The behavior in case $\alpha>\beta$ is different from the one in Theorem \ref{thm:main}. We claim that if $\alpha>\beta$ then, suitably normalized process $(N(Tt))_{t\ge 0}$ converges at least in the sense of finite dimensional distributions to a L\'evy process. We will discuss this case in a subsequent paper. 

The dychotomy between long range dependence for some parameters and independent increments for others is a known fenomenon in case of branching processes
 (see e.g. \cite{functlim3} and \cite{functlim4}). In the present case the intuition is that decreasing $\alpha$ means that some families will have descendants that are very far away from its ancestor, thus contributing to the dependence of increments in the limit. On the other hand, decreasing $\beta$ means that some particles may produce very many offspring, which increases the norming (note: if $\alpha<\beta$ we have $F_T=T^{\alpha+\frac {1+\alpha}{1+\beta}}$). Also, if $\beta$ is smaller, the families would have fewer generations until the time of extinction -- the tail of number of generations is lighter.
(Cf. \cite{AGH}, Proposition 2, where this was shown for a related branching mechanism \eqref{e:Gbeta}. Our branching in case of marked Hawkes process should behave similarly). 
Therefore it is natural to expect that if $\beta$ is sufficiently small in relation to $\alpha$, the dependence of increments will disappear in the limit.
\end{rem}

\medskip

The stable process $\zeta$ appearing in Theorem \ref{thm:main} has the following properties:
\begin{prop}\label{prop:continuity}
 The process $\zeta$ given by \eqref{e:limit_proc} has a modification with continuous sample paths and it is self-similar with index $H=\frac{1+\alpha(2+\beta)}{1+\beta}$.
\end{prop}

\medskip 
We will now discuss the functional convergence of the process $X_T$ given by \eqref{e:XTproc}. In what follows $D([0,\infty))$ is the Skorokhod space of c\`adl\`ag functions equipped with the usual $J_1$ topology.
We were not able to prove convergence $D([0,\infty))$  of the processes $X_T$ in the full generality of  Theorem \ref{thm:main}(b). However, we can do this under  an additional, technical assumption on $\varphi$, that seems to be quite natural. 

From \eqref{e:IR} it follows that under Assumption \ref{ass:phi_alpha}
 \begin{equation*}
\lim_{T\to \infty}  \frac 1{T^\alpha} \int_0^{tT}R(s)ds = c_\alpha t^\alpha
 \end{equation*} 
 Note that since $\alpha<1$ then $(t^\alpha-s^\alpha)\le  (t-s)^\alpha$ for $t\ge s$. It is therefore quite
 natural to expect that under some regularity condition on $\varphi$ we should have 
 \begin{equation}
\frac 1{T^\alpha}( I_R(Tt)-I_R(Ts))\le C (t-s)^{\alpha} \qquad 
  \text{for}\ t\ge s.
\label{e:estim_R}
  \end{equation}
Under a slightly weaker requirement we have convergence in the Skorokhod space.

\begin{thm}
\label{thm:tightness}
 Suppose that the assumptions of Theorem \ref{thm:main} are satisfied and additionally that for any $M>1$ there exist posistive constants $T_0,  C_M$ and $0<\varepsilon\le \alpha$  such that 
 \begin{equation}
 \frac 1{T^\alpha}(I_R(Tt)-I_R(Ts))\le C (t-s)^\varepsilon 
\qquad \text{for all}\ s,t \in [0,M], T\ge T_0.
   \label{e:assumption_tightness}
 \end{equation}
Then in part (b) of Theorem \ref{thm:main} we have convergence in law in $D([0,\infty))$ equipped with $J_1$ topology.
\end{thm}

\begin{ex}(a)
Suppose that $\varphi$ has the form
\begin{equation*}
\varphi(t)=\theta t^{\alpha-1}E_{\alpha, \alpha}(-\theta t^\alpha)
\end{equation*}
for some parameter $\theta>0$, where
 $E_{\alpha, \alpha}$ is the Mittag-Leffler function
\begin{equation*}
 E_{\alpha, \alpha}(t)=\sum_{k=0}^\infty\frac{t^k}{\Gamma(\alpha +\alpha k)}.
\end{equation*}
Then $R(s)=C s^{\alpha-1}$
(see Section 2.4.1 of \cite{HorstXu}) and consequently, \eqref{e:assumption_tightness} is satisfied.

\medskip
(b) Another simple example for which \eqref{e:assumption_tightness} holds is when $\varphi$ is the density of an $\alpha$-stable distribution totally skewed to the right. In this case we have $\varphi^{*k}(t)=k^{-\frac 1{\alpha}}\varphi(k^{-\frac 1\alpha}t)$, and $\varphi(s)\le  \frac {C_1}{(1+s^{\alpha+1})}$.
It is then not difficult to check that \eqref{e:estim_R} holds and therefore \eqref{e:assumption_tightness} is satisfied  with $\varepsilon=\alpha$.

\end{ex}

If $\nu$ has  finite variance, we obtain the following analogue of the result of Horst and Xu \cite{HorstXu}, recovering their result in case $\nu=\delta_1$.

\begin{thm}
 \label{thm:fin_variance}
Suppose $\varphi$ and $\nu$ satisfy the Assumptions \ref{ass:phi_alpha} and \ref{ass:nu_beta}($2$) and let 
\begin{equation}
F_T=T^{\frac{1+3\alpha} 2}.
 \label{e:FT2}
\end{equation}
a) $X_T$ defined by \eqref{e:XT} converge in law in
 $\SprimeR$ to a centered Gaussian $\SprimeR$-random variable $X$  with covariance functional
\begin{equation}
\label{e:cov_X}
E\<X,f_1\>\<X,f_2\>=\mu K\int_0^\infty G^{(\alpha)}f_1(x)G^{(\alpha)} f_2(x)x^\alpha dx
\qquad f_1, f_2\in \SR
\end{equation}
where
$$K=\frac 12 \int_0^\infty x^2\nu(dx)c_\alpha^{2+\beta}\alpha^{1+\beta}$$ with $c_\alpha$ given by \eqref{e:c_alpha}.\\
b) The processes $(X_T(t))_{t\ge 0}$ defined by \eqref{e:XTproc} converge in the sense of finite dimensional distributions to the process 
${\sqrt {2\mu K}}\zeta$ , where $\zeta$ is given by
\begin{equation}
\label{e:limit_procBM}
\zeta_t= \int_{[0,t]} u^{\frac \alpha{2}}(t-u)^\alpha B(du),
\end{equation}
where $B=(B_t)_{t\ge 0}$ is a standard Brownian motion. Moreover, if \eqref{e:assumption_tightness} is satisfied, then we have convergence in law in $D([0,\infty))$ equipped with $J_1$ topology.
\end{thm}

\begin{rem}\label{rem:BGT}
Returning to Theorem \ref{thm:main},
it is worth noting that the limiting $\SprimeR$-random variable $X$ in  Theorem \ref{thm:main}(b) 
has a form which is very similar to the limit process appearing in Theorem 2.1 (a) of \cite{functlim4} at a fixed time.
There, together with Bojdecki and Gorostiza we studied a different model of rescaled occupation time fluctuations of a branching particle system in which the particles moved according to symmetric $\alpha$-stable  L\'evy processes in $\R^d$, the lifetimes of particles were exponential and the branching mechanism had generating function
\begin{equation}
 G(s)=s+\frac{(1-s)^{1+\beta}}{1+\beta}
\label{e:Gbeta}
 \end{equation}
 with $0<\beta<1$,
which corresponds to a probability distribution  on $\Z_+$ satisfying \eqref{e:nu_beta}, therefore in the domain of attraction of $(1+\beta)$-stable law, totally skewed to the right. Even though the models are different and require different analysis the analogies are quite striking, especially when one takes $\varphi$ to be the $\alpha$-stable density, totally asymmetric.
If one changes somewhat the model of \cite{functlim4} by taking deterministic lifetimes (the branching occurs at  times $1,2,\ldots$), the offspring starts from the position of the parent and  the motion of particles is described by the $\alpha$-stable L\'evy process totally skewed to the right  and instead of considering the occupation time one looks at the empirical measure determined by the positions of all the parent particles at the branching event,  then the models become very similar.
The ``space'' in this model corresponds to ``time'' in the  marked Hawkes model studied in the present paper.

\medskip
In the branching model related to marked Hawkes processes the distance in time between the parent event and the child event has law with density $\varphi$, which is heavy-tailed. In this sense this model is also related to the age dependent branching particle systems considered in literature, for example in \cite{VatutinWakolbinger}, \cite{FVW} and \cite{LopezMimbela}. Howewer these articles considered problems different from the ones studied in the present paper.

\end{rem}
\bigskip

In the next theorem we show that the limitting behavior of the marked Hawkes process in Theorem \ref{thm:main} is the same as for an unmarked self-exciting process in which a single event induces a heavy tailed number of ``child'' events.
This corresponds to a branching mechanism which is not Poisson but which has generating function 
\eqref{e:Gbeta}. More precisely, we consider the following particle system in $\R$:

Suppose that the exogenous ($0$-generation) particles are distributed according to a Poisson random measure on $\Rp$ with intensity $\mu\ell$, that is the ``immigrants'' arrive according to a Poisson process with parameter $\ell$. 
Subsequently, each of the particles independently, gives rise to a random number of offspring, whose generating function is \eqref{e:Gbeta} with $0<\beta<1$.
The children particles are located to the right of the parent and their displacements in relation to the parent are given by i.i.d. random variables with density $\varphi$, $\norm{\varphi}_1=1$.
For $A\in \B(\R)$ we denote
\begin{equation*}
 N^\beta(A)=\{\# \text{of particles of all generations in set A}\}.
\end{equation*}
In the same way as before we have $N^\beta(t):=N^\beta([0,t])$, $X_T^\beta$ and $X_T^\beta(t)$ are defined accordingly, by formulas \eqref{e:XT} and \eqref{e:XTproc}. We have the following counterpart of Theorem \ref{thm:main}.

\begin{thm}
 \label{thm:beta_branching}
Suppose that $\varphi$ satisfies  Assumption \ref{ass:phi_alpha} with $\alpha<\beta$ and let $F_T$ be the norming given by \eqref{e:FT}.
Then \\
a) $X_T^\beta$ 
converge in law in $\SprimeR$ to an $\SprimeR$ valued $(1+\beta)$-stable random variable $X$ determined by
\eqref{e:f_char_X} with
$$K=\frac {1}{1+\beta}c_\alpha^{2+\beta}\alpha^{1+\beta}$$
b)  The processes $(X_T^\beta(t))_{t\ge 0}$ defined by \eqref{e:XTproc} converge in the sense of finite dimensional distributions to the process $\frac {(\mu K(-\cos \frac {\pi (1+\beta)}2))^{ 1/{1+\beta}}}\alpha \zeta$, where $\zeta$ is given by \eqref{e:limit_proc}. Moreover, if \eqref{e:assumption_tightness} holds, then there is convergence in law in $D([0,\infty))$ equipped with $J_1$ topology.
\end{thm}
From the proof it is evident that our method allows to extend further the class of possible branching mechanisms. See for example the remark at the beginning of the proof of Proposition \ref{prop:Lapl_conv}.

\section{Laplace transform and mean}\label{sec:Lapl_tr_and_mean}
In this section we derive a formula for the Laplace transform of the point measure $N$ associated to the times of events of  marked Hawkes process. For most of this section we do not need Assumptions \ref{ass:phi_alpha} and \ref{ass:nu_beta}, but only the basic assumptions stated at the beginning of the previous section. In particular, we assume that the system is critical. One can of course write similar, but slightly more complicated formulas in the non critical setting as well.
It is also easy to write the analogous formulas for the point measure $M$ involving both times of events and their marks, but we are not going do need them  here, since in our case the fact that the kernel function is of a multiplicative form \eqref{e:form_of_phi} simplifies the analysis.
These formulas are quite standard, using the branching property.

We will need some additional notation. 
Let $\rho$ denote the generic random variable describing  the number of offspring of a single particle (in one branching event). And let $G$ denote its generating function:
\begin{equation*}
 G(u)=E u^\rho \qquad \abs u\le 1.
\end{equation*}
Recall that we consider the system in which the mark of an event has distribution $\nu$ and that the conditional distribution of the number of offspring $\rho$, given the mark is Poisson with parameter equal to the mark. Also recall that the generating function of $Poiss(x)$ is $u\mapsto e^{x(1-u)}$, hence clearly
 \begin{equation}
  G(u)=\int_0^\infty e^{-x(1-u)}\nu(dx), \qquad \abs{u}\le 1.
  \label{e:generating_f}
 \end{equation}

Since we are only interested in the times of events in the marked Hawkes process, we can simplify the description of the empirical measure $N$ in the following way: we have ``immigrant'' particles distributed on $\Rp$ according to a Poisson random measure with intensity $\mu\ell$, each of the particles independently of each other has a random number of offspring with the generating function $G$ given by \eqref{e:generating_f}, the displacements of the offspring particles with relation to the parent are given by i.i.d. random variables with density $\varphi$. We stress that this simple form is due to the multiplicative form of the kernel function \ref{e:form_of_phi}.

Let us also denote
 \begin{align}
  \label{e:Hdef}
  H(u)=& G(1-u)-1+u \\
  =&\int_0^\infty \left(e^{-ux}-1+ux\right) \nu(dx).\label{e:HPoisson} \qquad 0\le u\le 1
 \end{align}

In the next proposition we describe the Laplace transform and the mean of the particle system $N$. In what follows $N^t$ stands for empirical measure of  an analogous branching particle system started by a single particle placed at $t$. The particle undergoes branching with displacements as described above. $N^t$ corresponds to the ``root'' particle and all of its descendants. 
\begin{prop}
 \label{prop:Laplace_2}
  Let  $f:\Rp\to \Rp$ be a Borel function. Then\\
a)
 \begin{equation}
  Ee^{-\<N, f\>}=\exp\left\{ -\int_0^\infty  g_{ f}(t)\mu dt\right\}
 \label{e:Laplace_2}
 \end{equation}
where
\begin{equation}
  g_{f}(t):=1-Ee^{-\<N^t, f\>}.
\label{e:def_g}
 \end{equation}
 Moreover, the function $g_f$
 satisfies the equation
\begin{align}
  g_{ f}(t)=&1-e^{-f(t)}
G\left(1-\int_0^\infty g_{f}(t+s)\varphi(s)ds\right)
 \label{e:g_eq}\\
 =& 1-e^{-f(t)}+e^{-f(t)}\int_0^\infty g_f(t+s)\varphi (s)ds -e^{-f(t)}H \left(\int_0^\infty g_f(t+s)\varphi(s)ds\right)
 \label{e:g_eqH}
\end{align}
b)
 \begin{equation}
  E\<N,f\>=
  \int_0^\infty  h_{ f}(t)\mu dt,
 \label{e:mean_2}
 \end{equation}
 where
 \begin{equation}
   h_{ f}(t):= E\<N^t,f\>
   \label{e:defh}
 \end{equation}
satisfies the equation
\begin{equation}
\label{e:renewal}
 h_{ f} (t) =  f(t)+\int_0^\infty  h_{f}(t+s)\varphi (s)ds.
\end{equation}
Moreover, $h_f$ may be expressed explicitly as
\begin{equation}
 h_f(t)=f(t)+\int_0^\infty f(t+s) R(s)ds,
 \label{e:hR}
\end{equation}
where $R$ is the resolvent defined in \eqref{e:R}.\\
c) For nonnegative $f$ the functions $g_f$ and $h_f$ satisfy
\begin{equation}
 0\le g_f(t)\le h_f(t)\qquad t\ge 0.
 \label{e:ineq_g_f}
\end{equation}
\end{prop}

We may also express \eqref{e:hR} in terms of $(\xi_j)_j$ -- i.i.d. random variables  with density $\varphi$ since we have
\begin{equation}
 \label{e:altR}
 \int_0^\infty f(t+s)R(s)ds=\sum_{k=1}^\infty E f(t+\sum_{j=1}^k \xi_j).
\end{equation}

\proofof{Proposition \ref{prop:Laplace_2}} 

Part a)
We can consider $N$ as a sum of families started by particles of generation $0$, that is 
\begin{equation*}
 N=\sum_{j=1}^\infty N^{t_j}_j
\end{equation*}
where $t_j, j=1,2,\ldots$ are the points of a Poisson random measure on $\Rp
$ with intensity $\mu \ell$ and $N^{t_j}_j$ are conditionally independent, given the initial Poisson random measure.

Therefore, to prove a) we condition on 
the initial Poisson random measure and use
the form of the Laplace transform of a Poisson random measure (see e.g. \cite{PeszatZabczyk}, Theorem 6.6)
\begin{align*}
 Ee^{-\<N, f\>}=&E\exp\left\{ e^{-\sum_j \<N^{t_j}_j,f\>}\right\}\\
 =& E\left\{ \prod_j E\left(e^{-\<N^{t_j},f\>}\big| (t_j)_j\right)\right\}\\
 =&\exp\left\{-\int_0^\infty(1-Ee^{-\<N^t, f\>})\mu dt\right\}.
\end{align*}
This gives \eqref{e:Laplace_2} with $g_{ f}$ defined by \eqref{e:def_g}.

To see that $g_{f}$ satisfies \eqref{e:g_eq} we 
use the branching property and condition on the first branching event. Suppose that the ``root'' particle is at $t$. This particle has $k$ offspring with probability $p_{k}$, where $\sum_{k=0}^\infty p_{k}s^k=G(s)$. The positions of the offspring are given by $t+\xi_j$, where $\xi_j$ are i.i.d. with density $\varphi$. Each of the offspring starts its own family that evolves in the same way, independently of all the other particles in the system. Therefore
\begin{align*}
 g_{ f}(t)=&1-e^{- f(t)}\sum_{k=0}^\infty p_{k}
 \left(Ee^{-\<N_j^{t+\xi}, f\>}\right)^k\\
 =& 1-e^{- f(t)}
G \left(Ee^{-\<N^{t+\xi}, f\>}\right)\\
=& 1-e^{- f(t)}
G \left(1- E g_{f}(t+\xi, \eta)\right).\\
\end{align*}
This proves \eqref{e:g_eq}.  \eqref{e:g_eqH} follows directly from \eqref{e:g_eq} and \eqref{e:Hdef}.
\medskip

 b) The proof of \eqref{e:mean_2}-\eqref{e:renewal}
is very similar but simpler  -- we skip it.

To prove that  \eqref{e:hR} is satisfied, assume first that $f$ is bounded and has compact support, $\supp f\subset[0,a]$, for some $0<a<\infty$.
Note that by  definiton of $h_f$ (see \eqref{e:defh}) we also have $\supp h_f\subset [0,a]$, since if $t>a$, then $N^{t}$ does not charge $[0,a]$. $h_f$ is bounded since $f$ is bounded and the mean number of particles that descend from a single particle  at $t$ and fall into $[0,a]$ is finite. The latter statement is completely clear if $r_a=\int_0^a\varphi(s)ds<1$ (as it is under Assumption \ref{ass:phi_alpha}), since then the expected number of particles of the $k$-th generation that fall into $[0,a]$ is estimated by $r_a^k$. Hence the total number of particles in $[0,a]$ is at most $ \sum_{k=0}^n r_a^n=\frac 1{1-r_a}<\infty$. 
Howewer, in this Proposition we do not assume \eqref{ass:phi_alpha}. In any case, there exists $b>0$ such that $r_b=\int_0^b\varphi(s)ds<1$ and $n$ such that $nb\ge a$. Then, by a similar argument, the expected total number of particles that fall into $[0,bn]$ is finite and bounded by $\sum_{k=1}^n (\frac 1{1-r_b})^k$. This proves that if $f$ is bounded and has bounded support then $h_f$ is bounded and has bounded support.

As usual, let $\xi_1,\xi_2, \ldots$  denote i.i.d. random variables with density $\varphi$. By  iteration of \eqref{e:renewal}
we have that for any $n$
\begin{equation*}
 h_f(t)=f(t)+\sum_{k=1}^n Ef(t+\xi_1+\ldots+\xi_n)+E h_f(t+\xi_1+\ldots+\xi_{n+1}).
\end{equation*}
By the strong law of large numbers $\sum_{k=1}^{n+1} \xi_k$ tends to infinity almost surely as $n\to \infty$, hence the last term in the above equation converges to zero by dominated convergence theorem.  The second last term converges to $\int_0^\infty f(t+s)R(s)ds$, therefore \eqref{e:hR} is satisfied for $f$ bounded with compact support. 

To prove \eqref{e:hR} in the general case of $f:\Rp\mapsto \Rp$ it is now sufficient to approximate $f$ by $f_n:=(f\wedge n)\ind_{[0,n]}$ and use monotone convergence.

\medskip

Part c) follows from the elementary inequality $1-e^{-x}\le x$ for $x\ge 0$ and the definitions of $g_f$ and $h_f$ given in \eqref{e:def_g} and \eqref{e:defh}, respectively. 
\qed

\bigskip

From now on we suppose that \textbf{Assumption \ref{ass:phi_alpha}} is satisfied. 

First we want to determine when \eqref{e:mean_2}--\eqref{e:hR} make sense for functions $f$ which are not necessarily nonnegative.
Recall \eqref{e:IR_def}, \eqref{e:c_alpha} and \eqref{e:IR}.
By Fubini theorem
\begin{align}
 \int_0^\infty \int_0^\infty\abs{f(t+s)}R(s)dsdt=&
 \int_0^\infty\int_s^\infty \abs{f(t)}dt R(s)ds\notag\\
 =&\int_0^\infty \abs{f(t)}I_R(t)dt.
 \label{e:Fubini}
\end{align}
Hence  
from  \eqref{e:hR} and \eqref{e:IR} we see that $E\<N,\abs{f}\><\infty$ for any $f\in \bF_\gamma$ with $\gamma > \alpha+1$.
Therefore also $E\<N,f\>$ is well defined. 

As a  direct consequence of  Proposition \ref{prop:Laplace_2} we obtain:

\begin{prop}
 \label{prop:Laplace_scaled}
Suppose that Assumption \ref{ass:phi_alpha} is satisfied and  let $f\in \bF_\gamma^+$ for $\gamma>\alpha+1$. Then
 \begin{equation}
  Ee^{-(\<N, f\>-E\<N,f\>)}=\exp\left\{ \int_0^\infty w_f(t)\mu dt\right\}
 \label{e:Laplace_centered}
 \end{equation}
where the function $w_f$ is defined by  
\begin{equation}w_f(t)=h_f(t)-g_f(t), \qquad t\ge 0
\label{e:def_wf}
\end{equation}
with $h_f$ and $g_f$ as in Proposition \ref{prop:Laplace_2}.
We have 
\begin{equation}
\label{e:ineq_wf}
 0\le w_f(t)\le h_f(t).
\end{equation}
Moreover, $w_f$  satisfies 
the equation
\begin{equation}
  w_f(t)=\int_0^\infty w_f(t+s)\varphi(s)ds + U_f(t),
 \label{e:eq_wf}
\end{equation}
where
\begin{equation}
U_f(t)= (f(t)- 1+e^{-f(t)})
  +(1-e^{-f(t)})\int_0^\infty g_f(t+s)\varphi(s)
ds
+e^{-f(t)}H
 \left(\int_0^\infty g_{f}(t+s)\varphi(s)ds\right),
 \label{e:Uf}
\end{equation}
with $H$ given by \eqref{e:Hdef}. 
Furthermore, 
$w_f$ may be also written as
\begin{equation}
 w_f(t)=U_f(t)+\int_0^tU_f(t+s)R(s)ds.
\label{e:UfR}
 \end{equation}
\end{prop}

\proof
From \eqref{e:Laplace_2} and \eqref{e:mean_2} and   we obtain
\eqref{e:Laplace_centered}. The equation \eqref{e:eq_wf} follows from  \eqref{e:g_eqH} and \eqref{e:renewal}.

Observe that $U_f$ is nonnegative, since for $x>0$ we have $x-1+e^{-x}=\int_0^x (1-e^{-u})du\ge 0$, hence also $H(x)$ is nonnegative. Therefore, \eqref{e:UfR} follows by  the standard iterative argument as in part b) of the proof Proposition \ref{prop:Laplace_2},
provided that we show that $w_f$ is bounded and vanishes at $+\infty$. By \eqref{e:ineq_wf}, it suffices to show that $h_f$ has these properties. Note that in the proof of Proposition \ref{prop:Laplace_2} we have already shown it in case when $f$ has compact support. In the general case of $f\in \bF_\gamma^+$ we can use \eqref{e:hR}. By \eqref{e:Fgamma} we have
 \begin{align*}
 \int_0^\infty f(t+s)R(s)\le  & 
C \frac 1{(1+t)^{\gamma-1-\alpha-\varepsilon}}
 \sum_{k=1}^\infty\int_{k-1}^{k} \frac 1{(1+s)^{1+\alpha+\varepsilon}}R(s)ds\\
 \le &C \frac 1{(1+t)^{\gamma-1-\alpha-\varepsilon}}
 \sum_{k=1}^\infty \frac 1{k^{1+\alpha+\varepsilon}}I_R(k)ds\\
 \le& \frac {C_1}{(1+t)^{\gamma-1-\alpha-\varepsilon}}.
\end{align*}
Where the last inequality follows from the fact that $I_R(k)\le C(1+k^\alpha)$, by \eqref{e:IR} and since $I_R$ is nondecreasing. 
By \eqref{e:hR} we see that $h_f$ is bounded and vanishes at $+\infty$, hence  $w_f$ also has these properties.
\qed

\section{Proofs of the main results}
\label{sec:proofs}
\subsection{Auxilliary lemmas}
\label{sec:auxilliary}
In the proof of the main result we will have to study the Laplace transform \eqref{e:Laplace_centered} with $f$ replaced by $f_T$ given by \eqref{e:fT} and we need to understand the asymptotics of various terms appearing in the Laplace transform. To prepare for this we formulate several lemmas.
The first one is an important consequence of the convergence \eqref{e:IR}:
\begin{lem}\label{lem:FR}
 Suppose that Assumption \ref{ass:phi_alpha} is satisfied and that $F\in \bF_\gamma$ for some $\gamma>1+\alpha$. 
 Then for any $0<\varepsilon<\gamma-1-\alpha$ there exists  $C=C(F,\varepsilon, \alpha)>0$ such that 
\begin{equation}
\label{e:Fbound}
 \sup_{T\ge 1}\int_0^\infty \abs{F(s+t)}T^{1-\alpha}R(sT)ds\le\frac C{(1+|t|)^{\gamma-1 -\alpha-\varepsilon}}.
\end{equation}
Moreover
\begin{equation}
 \lim_{T\to \infty}\int_0^\infty F(t+s)T^{1-\alpha} R(sT)ds=\alpha c_\alpha\int_0^\infty F(s+t)s^{\alpha-1}ds=\alpha c_\alpha G^{(\alpha)}F(t).
 \label{e:limF}
\end{equation}
\end{lem}

\begin{rem}
 \label{rem:heuristics}
Before we go to the formal proof let us observe that it is quite easy to understand why under Assumption  \ref{ass:phi_alpha} we should have convergence \eqref{e:limF}.

Recall that $\xi_1,\xi_2, \ldots$ denote a sequence of i.i.d. random variables with density $\varphi$.
Making a change of variables and then recalling the  the definition of $R$ (see \eqref{e:R} and \eqref{e:altR}) we can write
 \begin{align*}
  \int_0^\infty F(t+s)T^{1-\alpha}R(sT)ds
  =&\frac 1{T^\alpha}\int_0^\infty \sum_{k=0}^\infty F(t+
  \frac sT)R(s)ds\\
  &\sum_{k=0}^\infty
  \frac 1{T^\alpha}EF(t+\frac{\xi_1+\ldots+\xi_k }{T})\\
  =&\int_0^\infty \frac 1{T^\alpha}E F(t+\frac{\xi_1+\ldots+\xi_{\lceil x\rceil} }{T})dx\\
  =&\int_0^\infty E F(t+\frac{\xi_1+\ldots+\xi_{\lceil T^\alpha x\rceil} }{T})dx.
 \end{align*}
The processes $\frac {\xi_1 +\ldots +\xi_{\lceil T^\alpha x\rceil} }{T}$, $x\ge 0$ converge in law to the L\'evy process $C L_{\alpha}$, as $T\to \infty$ (see \cite{Feller} Ch. XVII.5, Theorem 2 on p. 577, and \cite{Skorokhod57}, Theorem 2.7). Therefore we expect that the last integral should converge to
\begin{equation*}
 \int_0^\infty E F(t+CL_{\alpha}(x))dx=C_1 G^{(\alpha)}F(t).
\end{equation*}
\end{rem}

\proofof{of Lemma \ref{lem:FR}}
First observe that from \eqref{e:Fgamma} it follows that for any $s,t\ge 0$ we have
\begin{equation}
 \label{e:3.1}
\abs{F(s+t)}\le \frac C{(1+t)^{\gamma-1-\alpha-\varepsilon}}\frac 1{(1+s+t)^{1+\alpha+\varepsilon}}.
 \end{equation}
Moreover
\begin{align}
 \int_0^\infty \frac 1{(1+s)^{1+\alpha+\varepsilon}}T^{1-\alpha}R(sT)ds\le &\sum_{k=1}^\infty\frac 1{k^{1+\alpha+\varepsilon}}\int_{k-1}^{k}T^{1-\alpha}R(sT)ds\notag\\
 \le &\sum_{k=1}^\infty\frac 1{k^{1+\alpha+\varepsilon}}
 T^{-\alpha} I_R(Tk).\label{e:3.2}
\end{align}
From \eqref{e:IR} it follows that there exists  $T_0$ such that 
\begin{equation*}
 \sup_{T\ge T_0} T^{-\alpha}I_R(T)\le 2c_\alpha,
\end{equation*}
therefore in \eqref{e:3.2},
if $Tk\ge  T_0$ then 
\begin{equation*}
  T^{-\alpha}I_R(Tk)\le 2c_\alpha k^\alpha.
\end{equation*} 
If  $Tk<T_0$ we estimate $T^{-\alpha}I_T(Tk)\le I_R(T_0)$
This shows that \eqref{e:3.2} is  bounded uniformly in $T\ge 1$ and \eqref{e:Fbound} follows from \eqref{e:3.2} and \eqref{e:3.1}.

\medskip
Now we proceed to $\eqref{e:limF}$. Some parts of this argument are similar to Remark 4.1 and Lemma 4.2 i n \cite{LopezMimbela}. From \eqref{e:Fgamma} it follows directly that for any $t\ge 0$, the function $\Rp\ni s\mapsto F(s+t)$ also belongs to $\bF_\gamma$, hence it is enough to prove the first equality in \eqref{e:limF}  for $t=0$.

Using estimates similar to \eqref{e:3.2} we have that 
\begin{equation*}
\sup_{T\ge T_0}\int_a^\infty \abs{F(s)}T^{1-\alpha}R(sT)ds\le \frac C{a^\varepsilon}\to 0 \qquad \textrm{as\ } a \to \infty.
\end{equation*}
Thus it suffices to prove \eqref{e:limF} for $F$ that has its support in a compact set $[0,a]$ for some $a>0$.

From \eqref{e:IR} it follows that for any $0\le t\le a$ we have
\begin{equation}
 \label{e:3.3}
 \lim_{T\to \infty}\frac{T^{-\alpha}I_R(tT)}{T^{-\alpha}I_R(aT)}
 =\left(\frac t a\right)^\alpha.
\end{equation}
Hence the probability measure with density 
\begin{equation*}
\frac{T^{1-\alpha}R(sT)}{T^{-\alpha}I_R(aT)}\ind_{[0,a]}(s)
\end{equation*}
(with respect to Lebesgue measure)
converges in law to the probability measure with density $\frac \alpha{a^\alpha}s^{\alpha-1}\ind_{[0,a]}(s)$, since \eqref{e:3.3} gives convergence of the corresponding distribution functions.
$F$ is bounded and has at most finite number of points of discontinuity and the limiting distribution does not charge these points, therefore a standard argument implies that
\begin{equation*}
 \lim_{T\to \infty}\frac {T^\alpha}{I_R(aT)}\int_0^\infty F(s)T^{1-\alpha}R(sT)ds=\frac \alpha{a^\alpha}\int_0^aF(s)s^{\alpha-1}ds.
\end{equation*}
Since $\lim_{T\to \infty}T^{-\alpha}I_R(aT)=c_\alpha a^\alpha$, this implies the first equality in \eqref{e:limF} for $t=0$.

The second equality in \eqref{e:limF} is just the definition of $G^{(\alpha)}$.
\qed

\medskip

As stated in Section \ref{sec:result}, our limit theorems correspond to the limits under rescaling \eqref{e:fT}. Lemma \ref{lem:FR} immediately leads to the following result on behavior of the function $h_{f_T}$ given by \eqref{e:defh} with $f_T$ given by \eqref{e:fT}. 

\begin{lem}  \label{lem:hTlim}
Suppose that Assumption \ref{ass:phi_alpha} is satisfied and that $f\in \bF_\gamma^+$ for some $\gamma>1+\alpha$ and let $f_T$ be given by \eqref{e:FT} with any positive norming $F_T$. Recalling \eqref{e:defh} let us denote $h_T:=h_{f_T}$. Then
 for any $\varepsilon<\gamma-1-\alpha$ there exists $C>0$ such that 
 \begin{equation}
  \label{e:hTbound}
\sup_{T\ge 1}  \frac {F_T}{T^\alpha}h_T(Tt)\le \frac C{(1+ t)^{\gamma-1-\alpha-\varepsilon}}\qquad t\ge 0.
 \end{equation}
Moreover, for any $t\ge 0$
\begin{equation}
\label{e:hTlim}
 \lim_{T\to \infty}\frac {F_T}{T^\alpha}h_T(Tt)=\alpha c_\alpha G^{(\alpha)}f(t).
\end{equation}
\end{lem}
Note that the form of the norming $F_T$ in this lemma is irrelevant. Also observe
that by \eqref{e:ineq_g_f}  the estimate \eqref{e:hTbound}  also gives the same estimates for the corresponding $g_{T}$ and $w_{T}$.

\proofof{Lemma \ref{lem:hTlim}}
From \eqref{e:hR} written for $f_T$ instead of $f$, we have that
\begin{align*}
 \frac{ F_T} {T^\alpha}h_T(Tt)=&T^{-\alpha} f(t)+
 T^{-\alpha}\int_0^\infty f(\frac{Tt+s}T) R(s)ds\\
 =&T^{-\alpha}f(t)+T^{1-\alpha}\int_0^\infty f(t+s)R(sT)ds.
\end{align*}
Hence both \eqref{e:hTbound} and \eqref{e:hTlim} follow from Lemma \ref{lem:FR}.
\qed

\bigskip
We also need to understand the behavior of $H$
 defined by \eqref{e:HPoisson} appearing in \eqref{e:eq_wf}
\begin{lem}
 \label{lem:H_properties_beta}
($\beta$) 
 Let $\nu$ satisfy Assumption \ref{ass:nu_beta}($\beta$) with $\beta\in(0,1)$ and let $H$ be given by \eqref{e:HPoisson}. Then  $H$ is nondecreasing and for any $x\ge 0$ we have
\begin{equation}
 \label{e:H_bound}
 H(x)= \frac{c_\nu}{\beta}\Gamma(1-\beta) x^{1+\beta}W(x),
\end{equation}
where $W$ is a nonnegative bounded function and satisfies
\begin{equation}
 \lim_{x\to 0+}W(x)=1.
 \label{e:limW}
\end{equation}
Moreover, there exists a constant $C>0$ such that for any $0\le x\le y$
\begin{equation}
 \label{e:H_bound_2}
 0\le H(y)-H(x)\le C(y-x)y^\beta.
\end{equation}

(2) Suppose that the probability law $\nu$ has mean $1$ and finite variance, then
 $H$ is nondecreasing and for any $x\ge 0$ we have
\begin{equation}
 \label{e:H_bound_var}
 H(x)= \frac{x^2} 2 W_2(x) \int_0^\infty u^2\nu(du),
\end{equation}
where $W_2$ is nonnegative, bounded by $1$  and satisfies
\begin{equation}
 \lim_{x\to 0+}W_2(x)=1.
 \label{e:limW_2_2}
\end{equation}
Moreover, there exists a constant $C>0$ such that for any $0\le x\le y$
\begin{equation}
 \label{e:H_bound_2_var}
 0\le H(y)-H(x)\le C(y-x)y.
\end{equation}
 \end{lem}
\proof
In both cases ($\beta$) and ($2$) we can write $H$ as
\begin{equation}
\label{e:H_alt}
 H(x)=\int_0^\infty \int_0^{ux}\int_0^w e^{-r}dr dw \nu(du). 
\end{equation}
From this we see directly that $H$ is nonnegative and nondecreasing in both cases. We still should prove that $H$ is finite, but this will be clear from the proof of \eqref{e:H_bound} and \eqref{e:H_bound_var}.

\medskip
Let us now consider the case ($\beta$).
We denote 
\begin{equation*}
 q_\nu(z)=z^{1+\beta}\int_z^\infty \nu(du).
\end{equation*}
From \eqref{e:nu_beta} it follows that 
\begin{equation}
 \sup_{z\ge 0} q_\nu (z)<\infty.
\label{e:sup}
 \end{equation}
This is clear, since for some $z_0>0$ from \eqref{e:nu_beta} it follows that for $z\ge z_0$ $q_\nu(z)\le 2c_\beta$, and  for $z\le z_0$ we have $q_\nu(z)\le z^{1+\beta}\le z_0^{1+\beta}$ since $\nu$ is a probability measure.

Hence using Fubini theorem in \eqref{e:H_alt} we obtain
\begin{align*}
 H(x)=&\int_0^\infty \int_r^\infty \int_{\frac wx}^\infty\nu(du) e^{ -r}dwdr\\
 =&x^{1+\beta}\int_0^\infty \int_r^\infty e^{-r}  \frac 1{{w}^{1+\beta}}q_\nu(\frac wx)dw dr.
\end{align*}
By Assumption \ref{ass:nu_beta},
 $q_\nu(\frac wx)$  converges to $c_\nu$ as $x\to \infty$ and since it is also bounded, both \eqref{e:H_bound} and \eqref{e:limW} follow from the fact that
\begin{equation*}
 \int_0^\infty \int_r^\infty e^{-r}  w^{-1-\beta} dw dr
 =\frac 1\beta\int_0^\infty r^{-\beta}e^{-r}dr =\frac 1{\beta}\Gamma(1-\beta).
\end{equation*}

Let us now consider \eqref{e:H_bound_2}. We have already shown the lower estimate. To see that the upper estimate is also satisfied we use \eqref{e:H_alt}, and to shorten the notation write $\eta$ for a random variable with law $\nu$
obtaining
\begin{align}
 H(y)-H(x)=&E\int_{\eta x}^{\eta y}\int_0^w e^{-r}dr dw
\le E\eta (y-x)\int_0^{\eta y}e^{-r}dr\label{e:4.13}\\
=&(x-y)\int_0^\infty e^{-r}E\left(\eta\ind_{\{\eta y\ge r\}}\right)dr.
\label{e:4.14}
 \end{align}
Now we look at the expectation under the integral in the last expression.  By Fubini theorem
\begin{equation}
 E\left(\eta\ind_{\{\eta y\ge r\}}\right)=\int_0^\infty P(\eta\ind_{\eta\ge \frac ry}>t)dt=\int_{\frac ry}^\infty P(\eta>t)dt=\int_{\frac ry}\frac 1{t^{1+\beta}}q_\nu(t)dt.
\label{e:4.15}
 \end{equation}
By \eqref{e:sup} this implies that
\begin{equation*}
  E\left(\eta\ind_{\{\eta y\ge r\}}\right)\le C_1\int_{\frac ry}^\infty \frac 1{t^{1+\beta}}dt\le \frac {C_1}\beta y^\beta r^{-\beta}.
\end{equation*}
Plugging this into \eqref{e:4.14} and using that $0<\beta<1$ we obtain \eqref{e:H_bound_2}.

\medskip
The case ($2$) is  simpler. We now assume that $\nu$ has a finite second moment. We again start with \eqref{e:H_alt} and make substitutions: $r'=\frac rx$, $w'=\frac wx$ obtaining that $H$ has a form
\begin{equation*}
 H(x)= x^2 
 \int_0^\infty \int_0^u\int_0^w e^{-rx}drdw \nu(du).
\end{equation*}
As $x\to 0$, the  integral above increases to $\frac 12\int_0^\infty u^2\nu(du)$, hence both \eqref{e:H_bound_var} and \eqref{e:limW_2_2} follow.

To obtain \eqref{e:H_bound_2_var} we use \eqref{e:4.13} and estimate $e^{-r}$ by $1$ in the last integral on the right hand side of \eqref{e:4.13}, using again finiteness of the second moment of $\nu$.
\qed

\subsection{Proof of Proposition \ref{prop:expectation}}

By \eqref{e:mean_2} and \eqref{e:renewal} the expectation $E\<N,f\>$ does not depend on $\nu$, hence it is the same as for $\nu=\delta_1$ corresponding to the not marked Hawkes process. Therefore,  Proposition \ref{prop:expectation} can be derived from the results of Horst and Xu \cite{HorstXu} 
(see Remark 2.3 and Proposition 2.6 therein).
However, since the proof follows directly from Lemma \ref{lem:FR} we will give this argument.

\medskip

From \eqref{e:mean_2} and \eqref{e:hR} we obtain that
\begin{align*}
\frac {1}{T^{1+\alpha}} E\<N, f(\frac \cdot T)\>
 =&\frac 1{T^{1+\alpha}}\int_0^\infty f(\frac tT)dt+ \frac 1{T^{1+\alpha}}\int_0^\infty \int_0^\infty f(\frac{t+s}{T}) R(s)dsdt\\
   =& \frac 1{T^{\alpha}}\int_0^\infty f(t)dt+  \int_0^\infty \int_0^\infty f(t+s) T^{1-\alpha}R(sT)dsdt.
\end{align*}
An application of Lemma \ref{lem:FR} finishes the proof.
\qed

\medskip

\subsection{Proof of Theorem \ref{thm:main}}\label{sec:proof}
\subsubsection{Convergence of Laplace transform and proof of Theorem \ref{thm:main}}\label{sec:proof_main}

In this section we reduce the problem to the following
convergence of the Laplace transforms of $\<X_T,f\>$ for nonnegative $f$. This is the key part of the proof of Theorem \ref{thm:main}. 

\begin{prop}
 \label{prop:Lapl_conv}
 Suppose that $f\in \bF^+_\gamma$.
 for some $\gamma>\frac {2+\alpha}{\beta(1+\beta)}+\frac{(2+\beta)(1+\alpha)}{1+\beta}$. 
 Then, under the assumptions of Theorem
 \ref{thm:main} we have
 \begin{equation}
  \label{e:Lapl_conv}
  \lim_{T\to \infty} E \exp(-\<X_T,f\>)=\exp\left\{\mu
  K\int_0^\infty \left(G^{(\alpha)}f(t)\right)^{1+\beta} t^\alpha dt
  \right\}
 \end{equation}
\end{prop}
\begin{rem}
 A straightforward calculation shows that for some constant $C>0$  the term in the exponent on the right hand side of \eqref{e:Lapl_conv} can be also written as
  \begin{equation*}
 C\int_0^\infty  G^{(\alpha)}\left[(G^{(\alpha)}f)^{1+\beta}\right](t)dt.
  \end{equation*}
\end{rem}

Proposition \ref{prop:Lapl_conv} will be proved in the next sections. Now we will show 
how this proposition implies Theorem \ref{thm:main}.

\proofof{Theorem \ref{thm:main}}
Let us recall that in Remark \ref{rem:main}b) we have observed that by \eqref{e:5.4a}, the process $\zeta$ may be written in the form \eqref{e:5.3}.
 $\zeta_t$ well defined since
\begin{equation*}
 \int_0^t\abs{G^{(\alpha)}\ind_{[0,t]}(y)}^{1+\beta}y^{\alpha}dy<\infty.
\end{equation*}
Moreover, formula (3.2.2) in \cite{ST} gives the finite dimensional distributions of $\zeta$. Namely, that for any $t_1,\ldots, t_n\ge 0$ and $\theta_1,\ldots, \theta_n\in \R$ 
\begin{equation*}
 E \exp \left(i\sum_{k=1}^n\frac {(K\mu)^{1/1+\beta}}{\alpha }\theta_k\zeta_{t_k}\right)
\end{equation*}
is equal to the right hand side of \eqref{e:f_char_X} with $f$ given by \eqref{e:f_sum}.
If additionally all $\theta_k$ are nonnegative then, the Laplace transform 
\begin{equation*}
 E\exp\left(-\sum_{k=1}^n\frac {(K\mu)^{1/1+\beta}}{\alpha }\theta_k\zeta_{t_k}\right)
\end{equation*}
is also expressed by \eqref{e:Lapl_X}.

Therefore, in both cases a) and b) the Laplace transform is given by the same formula \eqref{e:Lapl_X} with $f\in \bF^+$.

The remainder of the argument for the proof of Theorem \ref{thm:main} in case a) is exactly the same as in \cite{functlim3} (see Lemmas 3.4 and 3.5 therein, see also \cite{Iscoe}). 

In short, this argument is the following: First one proves convergence in law of $\<X_T,f\>$ for nonnegative $f\in \SR$, using convergence of Laplace transforms. In the general case, if $f$ is not necessarily positive, one writes $f$ as a difference of two nonnegative functions belonging to $\SR$ and  proves the convergence of the corresponding two dimensional distributions. This gives convergence of $\<X_T, f\>$ for $f\in \SR$ not necessarily positive. Moreover, this also shows  existence of an $\SprimeR$ random variable $X$ and convergence of $X_T$ to $X$ in $\SprimeR$, since the term in the exponent on the right hand side of \eqref{e:f_char_X} is continuous as a function form $\SR$
into $\C$ and by conuclearity of the space $\SprimeR$ (cf. \cite{Ito}). See \cite{functlim3} for more details.

In case  b), from convergence of Laplace transforms given by Proposition \ref{prop:Lapl_conv} with $f$ of the form  \eqref{e:f_sum} for all nonnegative $\theta_k$ we obtain convergence of the joint distributions of $(X_T(t_1),\ldots, X_T(t_k))$ to the stated distribution.

This finishes the proof of Theorem \ref{thm:main} provided that we show Proposition \ref{prop:Lapl_conv}. This will be done in the following section.
\qed

\subsubsection{Proof of Proposition \ref{prop:Lapl_conv}}.

 We  split the proof into several steps. Observe that concerning the function $H$ given by \eqref{e:Hdef} and \eqref{e:generating_f}, the proof only uses the properties stated in Lemma \ref{lem:H_properties_beta}, therefore we could replace $G$ by any other generating function such that the corresponding  $H$ has the same properties.

\medskip

\textbf{Step 1.} First we reformulate the problem by making use of the formulas developed in Section \ref{sec:Lapl_tr_and_mean}. Fix $f\in \bF^+_\gamma$, as in the assumptions. Recall \eqref{e:XT}, \eqref{e:fT}, \eqref{e:defh} and \eqref{e:def_g}.

To simplify the notation we  write
\begin{equation}\label{e:scaled}
 h_T(t):=h_{f_T}(t), \qquad g_T(t)= g_{f_T}(t), \qquad w_T(t)=h_T(t)-g_T(t), \qquad t\ge 0.
\end{equation}
Also recall that since $f$ is nonnegative, by \eqref {e:ineq_g_f} we have
\begin{equation}
0\le g_T(t)\le h_T(t), \qquad t\ge 0.
  \label{e:g_le_f}
\end{equation}

Rewriting Proposition \ref{prop:Laplace_scaled} in the new notation for the rescaled field $X_T$ we have:
 \begin{equation}
 E\exp\left(-\<X_T,f\>\right)=\exp\left(\mu\int_0^\infty w_T(t) dt\right),
 \label{e:Lapl_scaled}
\end{equation}
where 
\begin{equation}
  w_T(t)=U_T(t)+\int_0^\infty U_T(t+s)R(s)ds
  \label{e:w_TR}
\end{equation}
with
\begin{align}
  U_T(t)=& (f_T(t)- 1+e^{-f_T(t)})
  +(1-e^{-f_T(t)})\int_0^\infty g_T(t+s)\varphi(s)
ds\notag\\
&+e^{-f_T(t)}H
 \left(\int_0^\infty g_{T}(t+s)\varphi(s)ds\right),
 \label{e:UT}
\end{align}
where $H$ is defined by \eqref{e:Hdef}. 
Therefore, to prove \eqref{e:Lapl_conv} we need to show that
\begin{equation}
 \lim_{T\to \infty }
\int_{0}^\infty {w_T}(t)dt= K\int_0^\infty \left(G^{(\alpha)}f(t)\right)^{1+\beta} t^\alpha dt.
\label{e:lim_wT}
\end{equation}
The next steps will be devoted to showing \eqref{e:lim_wT}. The idea is the following: by changing variables we have
\begin{equation*}
 \int_0^\infty w_T(t)dt=\int_0^\infty Tw_T(Tt)dt.
\end{equation*}
We will show that $\sup_{T\ge 1}Tw_T(Tt)$ is bounded by an integrable function of $t$. Then we show that $Tw_T(Tt)$ behaves in the same way as $TS_T(Tt)$ where
\begin{equation*}S_T(t)=\int_0^\infty H(h_T(t+s))R(s)ds,
\end{equation*}
 and then prove that  $\int_0^\infty S_T(t)dt$ converges to \eqref{e:lim_wT}. In fact, we expect that 
$TS_T(Tt)$ also converges pointwise to $Const G^{(\alpha)}\left((G^{(\alpha)}f)^{1+\beta}\right)(t)$, but due to technical difficulties involving nested passages to the limit we were not able to prove formally this pointwise convergence.

\medskip

\textbf{Step 2.}
We will now show that 
for any $0<\kappa<\gamma (1+\beta)- {(1+\alpha)}{(2+\beta)}$ there exists 
 $C>0$ such that 
 \begin{equation}
  \label{e:wTbound}
  \sup_{T\ge 1}  T w_T(Tt)\le \frac C{(1+t)^\kappa}
\qquad t\ge 0.
 \end{equation}

To this end we estimate the terms of $U_T$ defined in \eqref{e:UT}.
We have the following estimates:
\begin{align}
0\le & f_T(t) -1+e^{-f_T(t)}\le\frac 12 f_T^2(t)\label{e:5.18}\\
\int_0^\infty &g_T(t+s)\varphi(s)ds\le \int_0^\infty h_T(t+s)\varphi(s)ds\le h_T(t).
\label{e:estim_gTint_hT}\\
0\le & (1-e^{-f_T(t)})\int_0^\infty g_T(t+s)\varphi (s)ds\le f_T(t)h_T(t),\label{e:5.19}
\end{align}
where in \eqref{e:estim_gTint_hT} we have used \eqref{e:g_le_f} and \eqref{e:hR}.

Moreover, by Lemma \ref{lem:H_properties_beta} and \eqref{e:estim_gTint_hT}
we also have that 
\begin{equation}
H
 \left(\int_0^\infty g_{T}(t+s)\varphi(s)ds\right)\le C_1 h_T^{1+\beta}(t), \qquad t\ge 0.
 \label{e:5.20}
\end{equation}
Since $f_T\le\frac {C_1}{F_T}\le C_2$, $f_T\le h_T$ and $0<\beta<1$ the above estimates imply that there exists $C>0$ such that for any $T\ge 1$
\begin{equation}
 U_T(t)\le C h_T^{1+\beta}(t)
 \label{e:b1}
\end{equation}
Fix $\varepsilon>0$ such that $0<\kappa<(\gamma-1-\alpha-\varepsilon)(1+\beta)-1-\alpha -\varepsilon$. Then by \eqref{e:hTbound} and \eqref{e:FT}
\begin{equation}
 Th_T^{1+\beta}(Tt)\le C\frac{T^{1+\alpha(1+\beta)}}{F_T^{1+\beta}}\frac 1{(1+t)^{(\gamma-1-\alpha-\varepsilon)(1+\beta)}}=C_{1}T^{-\alpha}\frac 1{(1+t)^{(\gamma-1-\alpha-\varepsilon)(1+\beta)}}.
 \label{e:b2}
\end{equation}
 Using the above, Lemma \ref{lem:FR} and \eqref{e:FT}
 we also have
\begin{align}
 T\int_0^\infty h_T^{1+\beta}(Tt+s)R(s)ds&=T^2\int_0^\infty h_T^{1+\beta}(T(t+s))R(Ts)ds\notag\\
 &\le C 
 T^{1-\alpha}
 \int_0^\infty \frac 1{{(1+t+s)}^{(\gamma-1-\alpha-\varepsilon)(1+\beta)}}R(sT)ds\notag\\
& \le 
 C_1 \frac 1{(1+t)^{(\gamma-1-\alpha-\varepsilon)(1+\beta)-1-\alpha -\varepsilon}}.
 \label{e:b3}
\end{align}
Combining \eqref{e:w_TR} and \eqref{e:b1}-\eqref{e:b3} we obtain \eqref{e:wTbound}.

\textbf{Step 3.} In this step we show that \eqref{e:wTbound} implies that $g_T$ behaves in the same way as $h_T$. 

Observe, that for $F_T$ given by \eqref{e:FT} for  we have
$\frac {F_T}{T^{\alpha}} = T^{\frac{\alpha-\beta}{1+\beta}} T$, hence from \eqref{e:wTbound} and \eqref{e:hTlim} we obtain directly that
 \begin{equation}
  \frac {F_T}{T^\alpha}w_T(Tt)\le T^{\frac{\alpha-\beta}{1+\beta}} \frac C{(1+t)^\kappa}
\qquad t\ge 0, T\ge 1.
 \label{e:wTtbound_2}
 \end{equation}
 This leads to the following corollary.
\begin{cor}
 Under the assumptions of Theorem \ref{thm:main} we have
\begin{equation}
\label{e:gTlim}
 \lim_{T\to \infty}\frac {F_T}{T^\alpha}g_T(Tt)=\lim_{T\to \infty}\frac {F_T}{T^\alpha}h_T(Tt)=\alpha c_\alpha G^{(\alpha)}f(t).
\end{equation}
\end{cor}

\medskip
Note that  the assumption $\alpha<\beta$ of Theorem \ref{thm:main} is crucially important here, to obtain that $g_T$ has the same behavior as $h_T$, as shown in \eqref{e:gTlim}. This is not the case if $\alpha>\beta$.

\medskip

It is now becoming clear that in the last term in \eqref{e:w_TR} we can replace $U_T(t+s)$ by $H(h_T(t+s))$ and this is the most significant term. The remaining ones will be negligeable. This is shown in the next step

\textbf{Step 4.}
 We will show that
  \begin{equation}
  \lim_{T\to \infty}\int_0^\infty w_T(t)dt=\lim_{T\to \infty}\int_0^\infty\int_0^\infty H\left(h_T(t+s)\right)R(s)dsdt,
  \label{e:lim_wT_lim_H_hT}
 \end{equation}
whenever any of the limits exists.

By \eqref{e:UT} and \eqref{e:5.18}-\eqref{e:5.20}
\begin{equation*}
\abs{ U_T(t)-H(\int_0^\infty g_T(t+s)\varphi(s)ds)}\le  C\left(f_T^2(t) +h_T(t)f_T(t) +f_T(t)h_T^{1+\beta}(t)\right) \le C_1 f_T(t)h_T(t).
\end{equation*}
Moreover, using  \eqref{e:H_bound_2},
 \eqref{e:renewal} and \eqref{e:scaled} we obtain
\begin{align*}
 \abs{H(h_T(t))- H(\int_0^\infty g_T(t+s)\varphi(s)ds)}&\le
 C \left(h_T(T)-\int_0^\infty g_T(t+s)\varphi (s)ds\right) h_T^\beta(t)\\
& \le C\left (f_T(t)+ \int_0^\infty w_T(t+s) \varphi(s) ds\right) h_T^\beta(t).
\end{align*}
Hence, 
\begin{align*}
 T\abs{U_T(Tt)-H(h_T(Tt))}
 \le C T\left( f_T(Tt) +\int_0^\infty w_T(T(t+\frac s T))\varphi (s)ds \right) h_T^\beta(Tt).
 \end{align*}
By \eqref{e:hTbound} and \eqref{e:wTtbound_2}
this implies
\begin{equation*}
  T\abs{U(Tt)- H(h_T(Tt))}\le C ({T^{-\alpha}}+T^{\frac{\alpha-\beta}{1+\beta}} )\frac {T^{1+\alpha(1+\beta)}}{F_T^{1+\beta}}\frac 1{(1+t)^{\kappa \beta}},
\end{equation*}
where $\kappa$ is as in \eqref{e:wTbound}. By the assumption on $\gamma$ we may choose $\kappa$ in such a way that $\kappa \beta> 2+\alpha+\varepsilon$  for some  $\varepsilon>0$.

Similarly as in \eqref{e:b2} and \eqref{e:b3}
we have
\begin{equation*}
 T\int_0^\infty\abs{U(Tt+s)- H(h_T(Tt+s))}R(s)ds 
 \le  C ({T^{-\alpha}}+T^{\frac{\alpha-\beta}{1+\beta}} )\frac 1{(1+t)^{\kappa\beta-1-\alpha-\varepsilon}}
\end{equation*}
This together with \eqref{e:b1} and \eqref{e:b2}
imply that
\begin{equation*}
\abs{ Tw_T(tT)-T\int_0^\infty H(h_T(Tt+s))R(s)ds}\le C({T^{-\alpha}}+T^{\frac{\alpha-\beta}{1+\beta}} ) \frac 1{(1+t)^{\kappa\beta-1-\alpha-\varepsilon}},
\end{equation*}
with  $\kappa\beta-1-\alpha-\varepsilon>1$. This in particular imply \eqref{e:lim_wT_lim_H_hT}, by a simple change of variables and since $\alpha>\beta$.

\textbf{Step 5.} 
Finally we will show that 
\begin{equation}
 \lim_{T\to\infty}\int_0^\infty \int_0^\infty H(h_T(t+s))R(s)dsdt=K\int_0^\infty \left(G^{(\alpha)}f(t)\right)^{1+\beta} t^\alpha dt.
 \label{e:lim_HT}
\end{equation}
This together with \eqref{e:lim_wT_lim_H_hT} will finish the proof of 
\eqref{e:lim_wT} and of the whole proposition.

In fact, we expect that 
$Tw_T(Tt)$ converges to $Const G^{(\alpha)}((G^{(\alpha)}f)^{1+\beta})(t)$, since $\frac T{F_T}h_T(Tt)\sim G^{(\alpha)}f(t)$, $\frac {F_T^{1+\beta}}{T^{\alpha(1+\beta)}}H(h_T(Tt))\sim C (\frac{F_T}{T^\alpha}h_T(Tt))^{1+\beta}(t)$, then another application of Lemma \ref{lem:FR} should give $Const G^{(\alpha)}((G^{(\alpha)}f)^{1+\beta})(t)$, but the nested dependence on $T$ makes it difficult to formally proove this pointwise convergence, so we only show \eqref{e:lim_HT}.

By Fubini
\begin{align*}
 \int_0^\infty \int_0^\infty H(h_T(t+s))R(s)dsdt=&\int_0^\infty 
 H(h_T(t))I_R(t)dt\\
 =&\int_0^\infty 
  \frac {H(h_T(Tt))}{h_T^{1+\beta}(Tt)}
 \left(\frac{F_T}{T^\alpha} h_{T}(Tt)\right)^{1+\beta}\frac{ T^{1+\alpha+\alpha(1+\beta)}}{F_T^{1+\beta}} \frac {I_R(Tt)}{T^\alpha}dt 
\end{align*}
Note that $\sup_{t\ge 0}h_T(Tt)\le C\frac{T^\alpha }{F_T}\to 0$ as $T\to \infty$.  Also, $T^{-\alpha}I_R(Tt)\le C(1+t^\alpha)$  by \eqref{e:IR} and since $I_R$ is increasing. Hence \eqref{e:hTlim}, \eqref{e:H_bound}, \eqref{e:limW}, \eqref{e:IR} and the dominated convergence theorem finishes the proof of \eqref{e:lim_HT}. This completes the proof of \eqref{e:lim_wT} and of  Proposition \ref{prop:Lapl_conv}.\qed

\subsection{Proof of Proposition \ref{prop:continuity}}
Self similarity means that for any $a>0$ the process $(\zeta_{at})_{t\ge 0}$ has the same finite dimensional distributions as $(a^H\zeta_t)_{t\ge 0}$, and this is clear from the way that the process $\zeta$ was obtained in Theorem \ref{thm:main}, but it can be also be easily checked by a direct calculation. 

To prove the existence of a continuous  modification we use the Kolmogorov continuity criterion. Namely, we will show that for any $M\ge 1$ there exists $C_M>0$ such that for any $0\le s<t\le M$ and $\lambda>0$ we have
\begin{equation}
 \P(\abs{\zeta_t-\zeta_s}>\lambda)
 \le \frac {C_M}{\lambda^{1+\beta}}(t-s)^{1+\alpha\beta}, 
\label{e:4.45}
 \end{equation}
which implies existence of a continuous modification.

Let $0\le s<t\le M$.
By \cite{Breiman}, Proposition 8.29, for any $\lambda>0$ we have
\begin{equation}
 \P(\abs{\zeta_t-\zeta_s}>\lambda)\le\int_0^1 (1-\Re \E e^{i\frac \theta\lambda(\zeta_t-\zeta_s)})d \theta.
\label{e:Breiman}
 \end{equation}

As observed at the beginning of the proof of Theorem \ref{thm:main} in Section \ref{sec:proof_main}, the characteristic function $\E e^{i\frac\theta\lambda (\zeta_t-\zeta_s)}$ has the same form as the right hand side of \eqref{e:f_char_X} with $f=\frac \theta\lambda\ind_{(s,t]}$ and therefore
\begin{equation*}
 1-\Re \E e^{i\frac \theta \lambda (\zeta_t-\zeta_s)}
 =1-e^{-C_1\int_0^\infty\abs{G^{(\alpha)}f(x)}^{1+\beta}x^\alpha dx}\cos \left(C_2 \int_0^\infty\abs{G^{(\alpha)}f(x)}^{1+\beta}x^\alpha dx\right).
\end{equation*}
Using  $1-e^{-u}\le u$ and $1-\cos(u)\le u$ for $u\ge 0$, we obtain that for $ \theta>0 $ and $\lambda>0$
\begin{equation}
  1-\Re \E e^{i\frac \theta \lambda (\zeta_t-\zeta_s)}\le C_3 \int_0^\infty\abs{G^{(\alpha)}f(x)}^{1+\beta}x^\alpha dx
  \le \frac{\theta^{1+\beta}}{\lambda^{1+\beta}}
\int_0^\infty\abs{G^{(\alpha)}\ind_{(s,t]}(x)}^{1+\beta}x^\alpha dx.
\label{e:4.47}
  \end{equation}
Therefore, by \eqref{e:Breiman} and \eqref{e:4.47}, to obtain \eqref{e:4.45} it suffices to show that
\begin{equation}
 \int_0^\infty\abs{G^{(\alpha)}\ind_{(s,t]}(x)}^{1+\beta}x^\alpha dx\le C_M(t-s)^{1+\alpha \beta}.
 \label{e:4.48}
\end{equation}
By  \eqref{e:5.4a} and the fact that for $0<\alpha\le 1$ and any $a,b\ge 0$ we have $(a+b)^\alpha\le a^\alpha+ b^\alpha$ we obtain 
\begin{equation*}
 G^{\alpha}\ind_{(s,t]}(x)=
 \frac 1\alpha \left((t-x)_+^\alpha -(s-x)^\alpha\right)
\le 
\frac 1 \alpha (t-s)^\alpha.
\end{equation*}
Hence for $s\le t\le M$
\begin{align*}
 \int_0^\infty \abs{G^{(\alpha)}\ind_{(s,t]}(x)}^{1+\beta}x^\alpha dx\le &\frac 1{\alpha^\beta}
(t-s)^{\alpha\beta}\int_0^\infty G^{(\alpha)}\ind_{(s,t]}(x) x^{\alpha} dx\\
\le  &\frac 1{\alpha^\beta}
(t-s)^{\alpha\beta}\int_0^\infty \int_x^\infty (y-x)^{\alpha-1}\ind_{(s,t]}(y)x^\alpha dy dx\\
=&\frac 1{\alpha^\beta}
(t-s)^{\alpha\beta}\int_s^t \int_0^y(y-x)^{\alpha-1} x^\alpha dxdy\\
=&\frac 1{\alpha^\beta}
(t-s)^{\alpha\beta}\int_s^t  y^{2\alpha} B(\alpha, \alpha+1)dy\\
\le&  M^{2\alpha}\frac 1{\alpha^\beta}  B(\alpha, \alpha+1)(t-s)^{1+\alpha\beta},
\end{align*}
where $B$ denotes the usual Beta function. \eqref{e:4.48} is proved. By \eqref{e:4.47} and \eqref{e:Breiman} this implies \eqref{e:4.45} and therefore finishes the proof of existence of continuous modification.\qed

\subsection{Proof of Theorem \ref{thm:tightness}}
In Theorem \ref{thm:main} we have already shown convergence in the sense of finite dimensional distributions. 

According to Theorems 13.5 and 16.7 in \cite{Billingsley}, to prove convergence in law in $D([0,\infty))$ it suffices to show that for any $M>0$ there exist positive $C_M$, $T_0$, $\kappa$, such that for any $a\le b\le d\le M$, $T\ge T_0$ and $\lambda>0$ we have
\begin{equation}
 \P(\abs{X_T(b)-X_T(a)}\wedge \abs{X_T(d)-X_T(b)}\ge \lambda)\le C_M \left(\frac 1{\lambda^{2+2\beta}}+\frac 1{\lambda^{1+\beta}}\right)(d-a)^{1+\kappa}.
 \label{e:tightness}
\end{equation}
(Note that  Theorem 13.5 in \cite{Billingsley}  has a slightly different form, with the same power good for $\lambda\ge 1$ and $\lambda\le 1$, but it is clear from the proof that only small $\lambda$ are important.
See also Theorems 7.8, 8.6, 8.8 and (8.40) in \cite{EthierKurtz}.) 

Recall, that $X_T(b)-X_T(a)=\<X_T,\ind_{(a,b]}\>$.

In the same way, as one proves  \eqref{e:Breiman} (cf.\  Proposition 8.29 in \cite{Breiman})  one can 
show that  there exists $C>0$ such that for any  real valued random variables $U,V$ 
\begin{align}
 \P(|U|\wedge |V|\ge \lambda)\le &C \int_{-1}^1\int_{-1}^1\E (1-e^{i\frac{Uu}{\lambda}})(1-e^{i\frac{Vv}{\lambda}})dudv \label{e:estim1}\\
 =&C \int_{-1}^1\int_{-1}^1\left( \E e^{i\frac {Uu+Vv}\lambda} -\E e^{i\frac {Uu}{\lambda}} -\E e^{i\frac {Vv}{\lambda}}+1\right)dudv
 \label{e:tail_estim}
\end{align}
Indeed,
the  estimate \eqref{e:estim1} is a simple consequence of 
\begin{equation*}
 \ind_{\abs{x}\ge \lambda}\le \frac 1{c_1}\left( 1-\frac {\sin(\frac x\lambda)}{\frac x \lambda}\right)=\frac 1{2c_1}\int_{-1}^1(1-e^{i\frac {xr}{\lambda}})dr,
\end{equation*}
where $c_1=\inf_{\abs{x}\ge 1}(1-\frac {\sin x}{x})$.

We want to use \eqref{e:tail_estim} for 
$U=\<X_T, \ind_{(a,b]}\>$ 
and $V=\<X_T, \ind_{(b,d]}\>$.
Similarly as at the begining of  the proof of Proposition \ref{prop:Laplace_2}, using the Poissonian structure, for any $f:\R_+\mapsto \R$ with compact support we have
\begin{equation}
 \E e^{i\<X_T,f\>}=\exp\left\{-\int_0^\infty\left(1-\E e^{i\<N^{t}, f_T\>} +i\E \<N^{t},f_T\> \right)\mu  dt\right\}
\label{e:4.44a}
 \end{equation}
and the real value of the term in the exponent is smaller or equal to $0$.

We want to expand the exponents in \eqref{e:tail_estim} and use the fact that we have already developed a good understanding of Laplace transforms.

Denote
\begin{equation}
 \label{e:ITf}
 J_T(f):= \abs{\E e^{i\<X_T,f\>}-1  +\int_0^\infty \left(1-\E e^{i\<N^{t}, f_T\>} +i\E \<N^{t}, f_T\> \right)\mu  dt}.
\end{equation}
Using \eqref{e:4.44a}, then
\begin{equation*}
 \abs{e^z-1-z}\le \frac{ \abs{z}^2}{2}, \qquad z\in \C, \Re z\le 0
\end{equation*}
and finally
\begin{equation*}
 \abs{e^{iy}-1-iy}\le (\frac 12 y^2)\wedge (2\abs{y}), \qquad y\in \R
\end{equation*}
we have
\begin{align*}
J_T(f)
 &\le \frac 12
\left(\int_0^\infty\int_0^\infty  \abs{1-\E e^{i\<M^{t}, f_T\>} +i\E \<N^{t}, f_T\> } \mu dt\right)^2
 \\
 & \le 2 \left( \int_0^\infty \E\left(\abs{\<N^{t},f_T\>}\wedge \abs{\<N^{t}, f_T\>}^2\right)\mu dt \right)^2\\
 & \le 2 \left(\int_0^\infty \int_0^\infty \E\left(\<N^{t},| f_T|\>\wedge \<N^{t},|f_T|\>^2\right)\mu dt \right)^2.
\end{align*}
$J_T(f)$ may be further estimated using 
\begin{equation*}
y\wedge y^2\le C (y+e^{-y}-1), \qquad y\ge 0.
\end{equation*}
This leads to
\begin{equation*}
 J_T(f)\le C\left(\int_0^\infty \E\left(\<N^{t},\abs{f}_T\>+e^{-\<N^{t},\abs{ f}_T\>}-1\right)\mu dt  \right)^2
\end{equation*}
Recalling the definitions of $h_f$, $g_f$ and $w_f$ the latter estimate is simply
\begin{equation}
 J_T(f)\le C_1 \left(\int_0^\infty w_{\abs{f}_T}(t)dt\right)^2
 \label{e:fin_estimate_IT}
\end{equation}
We now set $U=\<X_T, \ind_{(a,b]}\>$ 
and $V=\<X_T, \ind_{(b,d]}\>$ in \eqref{e:tail_estim} and use \eqref{e:ITf} and  \eqref{e:fin_estimate_IT} three times, obtaining
\begin{equation}
 \P(\abs{X_T(b)-X_T(a)}\wedge \abs{X_T(d)-X_T(b)}\ge \lambda)\le C\left(II_T+III_T\right)
 \label{e:tail_estimate2}
\end{equation}
where 
\begin{equation}
 II_T=
 \int_{-1}^1\int_{-1}^1 \int_0^\infty\E(1-e^{i\<N^{t},(\frac u\lambda \ind_{(a,b]})_T\>})(1-e^{i\<N^{t},(\frac v\lambda \ind_{(b,d]})_T\>}) dt dudv
 \label{e:IIT}
\end{equation}
and 
\begin{align}
 \abs{III_T}\le & C \int_{-1}^1\int_{-1}^1
\left(\int_0^\infty w_{\abs{\frac u\lambda \ind_{(a,b]} +\frac v\lambda \ind_{(b,d]}}_T}(t)dt\right)^2 dudv\notag\\
&+C \left(\int_{-1}^1\int_{-1}^1\left(\left(\int_0^\infty w_{\abs{\frac u\lambda \ind_{(a,b]}}_T}(t)dt\right)^2+ \left(\int_0^\infty w_{\abs{\frac v\lambda \ind_{(b,d]}}_T}(t)dt\right)^2\right)dudv\right).
\label{e:reszta}
\end{align}
Note that the terms  $i\E \<N^{t},(\frac u\lambda \ind_{(a,b]})_T\>$ and $i\E \<N^{t},(\frac v\lambda \ind_{(b,d]})_T\>$ in $II_T$ cancel out.

In \eqref{e:IIT} we use Fubini theorem and  integrate with respect to $u$ and $v$.
We use
\begin{equation*}
 \int_{-1}^1(1- e^{iux})du=\int_{-1}^1(1-\cos(ux))du\le 2\wedge x^2\le C (1-e^{-\abs{x}}) \qquad x\in \R,
\end{equation*}
then we apply
\begin{equation*}
 (1-e^{-x})(1-e^{-y})=(x+y)+e^{-(x+y)}-1 -(x+e^{-x}-1) -(y+e^{-y}+1),\qquad x,y \ge 0
\end{equation*}
and recall the definition of $w_f$,
obtaining
\begin{equation}
 0\le II_T=C\int_0^\infty \left(w_{(\frac 1{\lambda}(\ind_{(a,b]} +\ind_{(b,d]}))_T}(t)-
 w_{(\frac 1{\lambda}\ind_{(a,b]} )_T}(t)
 - w_{(\frac 1{\lambda}\ind_{(b,d]})_T}(t)
 \right)dt
\label{e:estim2}
 \end{equation}
From \eqref{e:tail_estimate2}, \eqref{e:reszta} and \eqref{e:estim2} it follows that to prove \eqref{e:tightness} it suffices to show that there exists $T_0\ge 1$, $C_M>0$ and $\varepsilon>0$ such that  for any 
nonegative function $f$ satisfying
$ f\le  \ind_{[a,d]}$ we have
\begin{equation}
 \int_0^\infty w_{(\frac{f}\lambda)_T}(t)dt \le \frac{C_M}{\lambda^{1+\beta}}(d-a)
 \label{e:w1}
\end{equation}
and for any  $0\le f_1\le \ind_{(a,b]}$ and $0\le f_2\le \ind_{(b,d]}$ we have
\begin{equation}
\label{e:w2}
 \int_0^\infty\left(w_{\frac 1\lambda (f_1+f_2)_T}(t)-w_{\frac 1\lambda (f_1)_T}(t)- w_{\frac 1\lambda (f_2)_T}(t)\right)dt\le {C_M}\left(\frac 1{\lambda^2}+\frac 1{\lambda^{1+\beta}}\right)(d-a)^{1+\varepsilon}.
\end{equation}
Let us start with \eqref{e:w1} and recall the notation $h_T=h_{f_T}$ and that $h_{(\frac 1{\lambda} f)_T}=\frac 1{\lambda} h_T$ by \eqref{e:hR}. From
 \eqref{e:b1}, \eqref{e:w_TR} 
and using Fubini we obtain
\begin{equation*}
\int_0^\infty w_{(\frac 1\lambda f)_T}(t)dt
\le C \frac 1{\lambda^{1+\beta}}\left(\int_0^\infty h_T^{1+\beta}(t)dt
+\int_0^\infty I_R(t) h_T^{1+\beta}(t)dt\right ).
\end{equation*}
The support of $f_T$ is contained in $[0,Td]\subset [0,TM]$, hence $h_T(t)=0$ for $t>TM$. Moreover,  $I_R$ is nondecreasing. Thus 
\begin{equation*}
 \int_0^\infty w_{(\frac 1\lambda f)_T}(t)dt
\le C \frac 1{\lambda^{1+\beta}} (1+I_R(TM))\int_0^\infty h_T^{1+\beta}(t)dt.
\end{equation*}
$M\ge 1$ is fixed, hence by \eqref{e:IR} we have that $I_R(TM)\le 2c_\alpha M^\alpha T^\alpha$ for sufficiently large $T$, $T\ge T_0$, therefore
\begin{equation}
 \int_0^\infty w_{(\frac 1\lambda f)_T}(t)dt
\le C \frac 1{\lambda^{1+\beta}}T^{\alpha}\int_0^\infty h_T^{1+\beta}(t)dt.
\label{e:wT3}
\end{equation}
Moreover, by \eqref{e:hR}
\begin{equation}
\frac{F_T}{ T^\alpha}  h_T(Tt)
 =\frac 1{T^{\alpha}}f(t)+\int_0^\infty f(t+s)T^{1-\alpha}R(sT)ds.
 \label{e:hT3}
\end{equation}
Since $\operatorname{supp} f\subset [0,M]$ this implies 
\begin{equation*}
 \frac{F_T}{T^\alpha}h_T(Tt)\le 1 +\frac {I_R(MT)}{T^\alpha}\le C M^\alpha  \qquad \text{for } t\ge 0, T\ge T_0
\end{equation*}
Recalling the definition of $F_T$ we have 
\begin{align}
 T^\alpha \int_0^\infty h_T^{1+\beta}(t)dt
 =&\frac{T^{1+\alpha}}{F_T^{1+\beta}}  
  T^{\alpha(1+\beta)}\int_0^\infty \left(\frac{F_T}{T^\alpha}h_T(tT)\right)^{1+\beta}dt\notag\\
 \le & C_1 M^{\alpha \beta} \int_0^\infty 
 \frac{F_T}{T^\alpha }h_T(Tt)dt\notag\\
 =&C_1 M^{\alpha\beta}\left(\frac 1{T^\alpha } (d-a)
 +\int_0^M f(t) \frac {I_R(Tt)}{T^\alpha}dt
 \right)\notag \\
 \le & C_1 M^{\alpha \beta}(\frac 1{T_0^\alpha}+CM^\alpha) (d-a) \label{e:estim_h}
\end{align}
Putting together \eqref{e:wT3} and \eqref{e:estim_h} we obtain \eqref{e:w1}. Note that in this part we have not used the additional assumption \eqref{e:assumption_tightness}. It will only be used in the proof of \eqref{e:w2}.

\medskip
We now turn to \eqref{e:w2}. 
Denote 
\begin{equation}
v_T(t)=w_{\frac 1\lambda (f_1+f_2)_T}(t)-w_{\frac 1\lambda (f_1)_T}(t)- w_{\frac 1\lambda (f_2)_T}(t)
 \label{e:vT}
\end{equation}
We use \eqref{e:w_TR} for the corresponding functions  and observe that there are some cancellations due to the fact that the supports of $f_1$ and $f_2$ are disjoint. We obtain
\begin{equation}
 v_T(t)=Z_T(t)+\int_0^\infty Z_T(t+s)R(s)ds,
\label{e:ZT}
 \end{equation}
 where
\begin{align*}
 Z_T(t)=&(1-e^{-(\frac {f_1}\lambda)_T(t)})\int_0^\infty \left(g_{(\frac {f_1+f_2}\lambda)_T}(t+s) - g_{(\frac {f_1}\lambda)_T}(t+s)\right)\varphi(s)ds\\
&+(1-e^{-(\frac {f_2}\lambda)_T(t)})\int_0^\infty \left(g_{(\frac {f_1+f_2}\lambda)_T}(t+s) - g_{(\frac {f_2}\lambda)_T}(t+s)\right)\varphi(s)ds\\
&+e^{-(\frac {f_1+f_2}\lambda)_T(t)}H(\int_0^\infty 
g_{(\frac {f_1+f_2}\lambda)_T}(t+s)\varphi(s)ds)\\
&-e^{-(\frac {f_1}\lambda)_T(t)}H(\int_0^\infty 
g_{(\frac {f_1}\lambda)_T}(t+s)\varphi(s)ds)\\
&-e^{-(\frac {f_2}\lambda)_T(t)}H(\int_0^\infty 
g_{(\frac {f_2}\lambda)_T}(t+s)\varphi(s)ds).
\end{align*}
$Z_T(t)$  vanishes for $t\ge dT$ and $I_R(MT)\le CM^\alpha T^\alpha$, hence by Fubini
\begin{equation}
 \int_0^\infty v_T(t)dt\le CM^\alpha T^\alpha \int_0^M \abs{Z_T(t)}dt.
\label{e:vTestim}
 \end{equation}
We will
estimate the absolute value of $Z_T$. The last three terms of $Z_T$ may be bounded  using that for $0\le f\le \ind_{(a,d]}$ by \eqref{e:H_bound}, \eqref{e:limW} and \eqref{e:g_le_f} we have
\begin{align*}
 H(\int_0^\infty g_{(\frac f \lambda)_T}(t+s)\varphi(s)ds)\le &C\left(\int_0^\infty h_{(\frac f \lambda)_T}(t+s)\varphi (s)ds\right)^{1+\beta}\\
 =&\frac C{\lambda^{1+\beta}} \left(\int_0^\infty h_{(\ind_{(a,d]})_T}(t+s)\varphi (s)ds\right)^{1+\beta}.
\end{align*}
In the first  term of $Z_T$ we use   the definition of $g_f$ (see \eqref{e:def_g}) and then \eqref{e:g_le_f} to obtain
\begin{equation*}
0\le g_{(\frac {f_1+f_2}\lambda)_T} -  g_{(\frac {f_1}\lambda)_T}\le g_{(\frac {f_2}\lambda)_T}
 \le
 h_{(\frac {f_2}\lambda)_T}\le \frac 1{\lambda}
 h_{(f_2)_T}.
\end{equation*}
Moreover, by  \eqref{e:renewal} and by \eqref{e:hR}  
\begin{equation*}
\int_0^\infty h_f(t+s)\varphi (s)ds=h_f(t)-f(t)=\int_0^\infty f(t+s)R(s)ds          
\end{equation*}
Hence
\begin{equation*}
 0\le \int_0^\infty \left(g_{(\frac {f_1+f_2}\lambda)_T}(t+s) - g_{(\frac {f_1}\lambda)_T}(t+s)\right)\varphi(s)ds\le \frac C{\lambda}\int_0^\infty 
 (f_2)_T(t+s)R(s)ds.
\end{equation*}
Using also 
\begin{equation*}
 (1-e^{-(\frac {f_1}\lambda)_T(t)})\le \frac 1{\lambda} (\ind_{(a,d]})_T(t)
\end{equation*}
and analogous estimates  for the second term of $Z_T$ we obtain
the following estimate on $\abs{Z_T}$:
\begin{equation}
 \abs{Z_T(t)}\le \frac C{\lambda^2}
  (\ind_{(a,d]})_T(t)  \int_0^\infty {(\ind_{(a,d]})_T}(t+s)R (s)ds
  +\frac C{\lambda^{1+\beta}} \left(\int_0^\infty {(\ind_{(a,d]})_T}(t+s)R (s)ds\right)^{1+\beta}.
  \label{e:ZT_estimate}
\end{equation}
For any $t\ge 0$ we have that 
\begin{align}
A_T:=\int_0^\infty (\ind_{(a,d]})_T({t+s})R (s)ds
 =&\frac {1}{F_T}\int_0^\infty\ind_{(Ta,Td]}( {t+s})R(s)ds\notag\\
=& \frac {1}{F_T}\left( I_R((Td-t)\vee 0)- I_R((Ta-t)\vee 0)\right)\notag\\
\le& \frac {T^\alpha}{F_T}(d-a)^\varepsilon.
\label{e:Jestim}
 \end{align}
by assumption \eqref{e:tightness}. Recall that we may assume that  $0<\varepsilon<1$. 
Plugging this into \eqref{e:ZT_estimate}
we obtain
\begin{equation}
 \abs{Z_T(t)}\le \frac C{\lambda^2}\frac{T^\alpha}{F_T^2}
  \ind_{(Ta,Td]}(t)  (d-a)^\varepsilon
  +\frac C{\lambda^{1+\beta}} 
  \frac{T^{\alpha\beta}}{F_T^{\beta}}(d-a)^{\varepsilon\beta} \int_0^\infty {(\ind_{(a,d]})_T}(t+s)R (s)ds
  \label{e:ZT_estimate2}
\end{equation}
By Fubini and \eqref{e:vT}-\eqref{e:ZT_estimate} we obtain
\begin{align*}
 \int_0^\infty v_T(t)dt
\le & \int_0^T\abs{Z(t)}(1+I_R(t))dt\\
\le & C(1+M^\alpha)\int_0^\infty \abs{Z(t)}dt\\
  \le & \frac{C_M T^{1+2\alpha}}{\lambda^2 F_T^2}(d-a)^{1+\varepsilon}
 +\frac {C_M}{\lambda^{1+\beta}} 
  \frac{T^{\alpha+\alpha\beta}}{F_T^{1+\beta}}(d-a)^{\varepsilon\beta}\int_0^\infty \ind_{(aT, dT]}(t)I_R(t)dt\\
  \le& \frac{C_M T^{1+2\alpha}}{\lambda^2 F_T^2}(d-a)^{1+\varepsilon}
 +\frac {C_M}{\lambda^{1+\beta}} 
  \frac{T^{1+2\alpha+\alpha\beta}}{F_T^{1+\beta}}(d-a)^{1+\varepsilon\beta}
\end{align*}
 Recalling the definition of $F_T$ finishes the proof of \eqref{e:w2} and the proof of tightness.
Note that the only point where we used the additional estimate \eqref{e:assumption_tightness} was to estimate $A_T$ in \eqref{e:Jestim}.\qed

\subsection{Proof of Theorem \ref{thm:fin_variance}}
The proof proceeds in the same way as the proof of 
Theorem \ref{thm:main} and \ref{thm:tightness} if we set $\beta=1$. Using the form of the Laplace transform of Gaussian random variables we reduce the proof to the proof of Proposition \ref{prop:Lapl_conv} with $\beta=1$. The argument for the latter one is the same, except that instead of  Lemma \ref{lem:H_properties_beta}$(\beta)$ we employ part  $(2)$ of the same lemma. We omit the details.\qed

\subsection{Proof of Theorem \ref{thm:beta_branching}}
The proof is essentially the same as that of Theorem \ref{thm:main}.  The generating function  $G$ given by \eqref{e:generating_f} is now replaced by  \eqref{e:Gbeta}.  We have the same formula for Laplace transform of  $X^\beta_T$ as in \eqref{e:Lapl_scaled} and \eqref{e:w_TR}, but with different $H$. Now
\begin{equation*}
 H(s)=G(1-s)-(1-s)=\frac{s^{1+\beta}}{1+\beta}.
\end{equation*}
We trivially have \eqref{e:H_bound} and \eqref{e:limW} with $\frac{c_\nu}{\beta}\Gamma(1-\beta)$ replaced by $\frac 1{1+\beta}$. Moreover \eqref{e:H_bound_2} holds with $C=1+\beta$.

The rest of the argument is the same (see the remark at the beginning of the proof of Proposition \ref{prop:Lapl_conv}). The argument for tightness is also the same. \qed

\bibliographystyle{plain}

\begin{thebibliography}{22}
\itemsep=0pt

\bibitem{AGH} J. Avan, N. Grosjean, T. Huillet,
On extreme events for non-spatial and spatial branching Brownian motions.
{\em Physica D: Nonlinear Phenomena},
Volumes 298–299 (2015),
 13--20.
 
\bibitem{BacryDelattreHoffmann} E.~ Bacry, S.~Delattre, M.~Hoffmann.
Some limit theorems for Hawkes processes and application to financial statistics.
{\em Stochastic Process. Appl.} 123 (2013), no. 7, 2475–2499.

\bibitem{Breiman}  L.~Breiman,  Probability. Addison–Wesley, Reading, MA., 1968.

\bibitem{Billingsley}
P.~P. Billingsley, Convergence of probability measures, second edition, 
Wiley Series in Probability and Statistics: Probability and Statistics A Wiley-Interscience Publication,  Wiley, New York, 1999.

\bibitem{Bingham}
N.H. Bingham, C.M. Goldie, J.L. Teugels, 
Regular variation.
Encyclopedia Math. Appl., 27
Cambridge University Press, Cambridge, 1989. xx+494 pp.

\bibitem{functlim3}T.~Bojdecki, L.~G. Gorostiza, and A.~Talarczyk.
A long range dependence stable process and an infinite variance branching system.
{\em Ann. Probab.} 35 (2007), no. 2, 500–527.
\bibitem{functlim4}
T.~Bojdecki, L.~G. Gorostiza, and A.~Talarczyk.
\newblock Occupation time fluctuations of an infinite-variance branching system in large dimensions. {\em
Bernoulli} 13 (2007), no. 1, 20–39.

\bibitem{DaleyVereJones}
D.J.~Daley, D.~Vere-Jones.
An introduction to the theory of point processes. Vol. I.
Elementary theory and methods. Second edition
Probab. Appl. (N. Y.)
Springer-Verlag, New York, 2003.

\bibitem{DGW} D.A.~Dawson, L.G.~Gorostiza, A.~ Wakolbinger.
Degrees of transience and recurrence and hierarchical random walks.
Potential Anal. 22 (2005), no. 4, 305–350.

\bibitem{EthierKurtz}
S.~N. Ethier and T.~G. Kurtz,  Markov processes, Wiley Series in Probability and Mathematical Statistics: Probability and Mathematical Statistics, Wiley, New York, 1986.

\bibitem{Feller}W. Feller, \textit{An Introduction to Probability Theory and Its Applications}, vol. 2, 2nd edn., Wiley,
New York, 1971.

\bibitem{FVW} K. Fleischmann, V.A. Vatutin, A. Wakolbinger, Branching systems with long-living particles at the critical dimension. Theory Probab. Appl., {\bf 47} (2023), no. 3, 429--454. 

\bibitem{Hardiman} S.J. Hardiman, N. Bercot, J.P. Bouchaud,\ Critical reflexivity in financial markets: a Hawkes process analysis. Eur. Phys. J. B  (2013) 86, 442. 



\bibitem{Hawkes1} A.G.~Hawkes.
Point spectra of some mutually exciting point processes.
{\em J. Roy. Statist. Soc.} Ser. B 33, 438–443, 1971.
\bibitem{Hawkes2} A.G.~Hawkes.
Spectra of some self-exciting and mutually exciting point processes.
{\em Biometrika} 58 (1971), 83–90.

\bibitem{HawkesRev}
A.~G.~Hawkes. Hawkes processes and their applications to finance: a review, Quant. Finance {\bf 18} (2018), no.~2, 193--198; 

\bibitem{HawkesOakes} 
A.G.~Hawkes, D.~Oakes.
A cluster process representation of a self-exciting process.
{\em J. Appl. Probability} 11 (1974), 493–503.

\bibitem{HorstXu_marked}
U.~Horst,W.~Xu.
Functional limit theorems for marked Hawkes point measures.
Stochastic Process. Appl. 134 (2021), 94–131.

\bibitem{HorstXu} U.~Horst, W.~Xu. Functional limit theorems for Hawkes processes. { \em Probab. Theory Relat. Fields} {\bf 194} (2026), 917–996.


\bibitem{HorstXuZhang} U.~Horst, W.~Xu and R.~Zhang, Convergence of heavy-tailed Hawkes processes and the microstructure of rough volatility," arXiv:2312.08784,  2023.

\bibitem{Iscoe} I.~Iscoe. A weighted occupation time for a class of measure-valued branching processes.
{\em Probab. Theory Relat. Fields} 71 (1986), no. 1, 85–116.
\bibitem{Ito} K.~It\^o.
Foundations of stochastic differential equations in infinite-dimensional spaces.
CBMS-NSF Regional Conf. Ser. in Appl. Math., 47, SIAM, Philadelphia,  1984. 

\bibitem{JaissonRosenbaum_nearly_unstable} T. Jaisson and M. Rosenbaum, Limit theorems for nearly unstable Hawkes processes, {\em Ann. Appl. Probab.} {\bf 25} (2015), no.~2, 600--631.

\bibitem{JR_2}T. Jaisson and M. Rosenbaum, Rough fractional diffusions as scaling limits of nearly unstable heavy tailed Hawkes processes, {\em Ann. Appl. Probab.} {\bf 26} (2016), no.~5, 2860--2882.


\bibitem{Laub}
P.~J. Laub, Y. Lee, P.K. Pollett, T. Taimre, Hawkes models and their applications, {\em Annu. Rev. Stat. Appl.} {\bf 12} (2025), 233--258. 

\bibitem{Li_Pang}
B. Li, G. Pang,
Functional limit theorems for nonstationary marked Hawkes processes in the high intensity regime.
{\em Stochastic Process. Appl. } {\bf 143} (2022), 285–339.

\bibitem{LopezMimbela}
J.~A. L\'opez-Mimbela, A. Murillo-Salas and J.~H. Ram\'irez~Gonz\'alez, Occupation time fluctuations of an age-dependent critical binary branching particle system, ALEA Lat. Am. J. Probab. Math. Stat. {\bf 21} (2024), no.~1, 593--625.


\bibitem{Ogata} Y. Ogata,  Statistical Models for Earthquake Occurrences and Residual Analysis for Point Processes, {\em
Journal of the American Statistical Association,} Vol. 83, No. 401 (Mar., 1988),  9--27.

\bibitem{PeszatZabczyk} S.~Peszat, J.~Zabczyk,
Stochastic partial differential equations with Lévy noise.
An evolution equation approach
Encyclopedia Math. Appl., 113
Cambridge University Press, Cambridge, 2007.

\bibitem{ST} Samorodnitsky, G. and Taqqu, M.S.,  {\it Stable Non-Gaussian Random Processes,} Chapman \& Hall, New York, 1994.


\bibitem{Skorokhod57} A.V. Skorokhod, Limit theorems for stochastic processes with independent increments, {\em Theory Probab. Appl. }, Vol.2, No.2 (1957), 138-171.

\bibitem{VatutinWakolbinger} V. Vatutin, A. Wakolbinger, Spatial Branching Populations with Long Individual Lifetimes.
{\em Theory Probab. Appl.}, { \bf 43} (4) (1999), 620–632. 
\end{thebibliography}

\end{document}